\documentclass{article}

\usepackage{PRIMEarxiv}

\usepackage[utf8]{inputenc} 
\usepackage[T1]{fontenc}    
\usepackage{hyperref}       
\usepackage{booktabs}       
\usepackage{nicefrac}       
\usepackage{microtype}      
\usepackage{fancyhdr}       
\usepackage{graphicx}       
\graphicspath{{media/}}     
 
\usepackage{amssymb}

\usepackage{natbib}
\bibpunct[, ]{(}{)}{;}{a}{}{,}

\usepackage{caption}
\usepackage{subcaption}
\usepackage{amsfonts}
\usepackage{amsmath}
\usepackage{tikz}
\usetikzlibrary{decorations.pathreplacing}
\usepackage{makecell, multirow}
\usepackage{rotating}
\usepackage{units}
\usepackage{url}
\usepackage{placeins}
\usepackage{float}
\usepackage{algorithm}
\usepackage{algpseudocode}

\newcommand{\Js}{\sum_{j\in\mathcal{J}}}
\newcommand{\Is}{\sum_{i\in\mathcal{I}}}

\newcommand{\Ss}{\sum_{s\in\mathcal{S}}}

\newcommand{\Ks}{\sum_{k\in\mathcal{K}}}
\newcommand{\ForJ}{\forall j \in \mathcal{J}}
\newcommand{\ForI}{\forall i \in \mathcal{I}}
\newcommand{\ForS}{\forall s \in \mathcal{S}}
\newcommand{\ForM}{\forall m \in \mathcal{M}}
\newcommand{\ForK}{\forall k \in \mathcal{K}}

\newcommand{\ForN}{\forall n \in \mathcal{N}}
\newcommand{\st}{\mathrm{s.t.}}
\usepackage{xspace}
\newcommand{\MP}{\texttt{MP\_<n>}\xspace}
\newcommand{\MPd}{\texttt{MP\_2}\xspace}
\newcommand{\MPt}{\texttt{MP\_3}\xspace}
\newcommand{\MPq}{\texttt{MP\_4}\xspace}
\newcommand{\MSt}{\texttt{MS3}\xspace}
\newcommand{\MStp}{\texttt{MS3\_<n>}\xspace}
\newcommand{\FOSVA}{\texttt{FOSVA}\xspace}
\newcommand{\TS}{\texttt{TS}\xspace}
\newcommand{\TSn}{\texttt{TS\_noS}\xspace}
\newcommand{\MStt}{\texttt{MS3\_3}\xspace}
\newcommand{\MStq}{\texttt{MS3\_4}\xspace}
\newcommand{\SSq}[1]{\texttt{SSt\_#1}\xspace}

\newcommand{\bc}{\begin{center}}
\newcommand{\ec}{\end{center}}

\title{Rolling horizon policies for multi-stage stochastic assemble-to-order problems
}

\author{
  Daniele Giovanni Gioia\thanks{Corresponding author},\ \ Edoardo Fadda, Paolo Brandimarte \\
  Department of Mathematical Sciences "Giuseppe Luigi Lagrange",\\ DISMA, Politecnico di Torino, Corso Duca degli Abruzzi 24, \\
  Turin (TO), Italy.\\
  \texttt{\{daniele.gioia, edoardo.fadda, paolo.brandimarte\}@polito.it}\\
  ORCID: 0000-0001-8979-4174, 0000-0002-5599-6349, 0000-0002-6533-3055.
}

\begin{document}
\maketitle

\begin{abstract}
Assemble-to-order approaches deal with randomness in demand for end items by producing components under uncertainty, but assembling them only after demand is observed. Such planning problems can be tackled by stochastic programming, but true multistage models are computationally challenging and only a few studies apply them to production planning. Solutions based on two-stage models are often short-sighted and unable to effectively deal with non-stationary demand. A further complication may be the scarcity of available data, especially in the case of correlated and seasonal demand. 
In this paper, we compare different scenario tree structures. In particular, we enrich a two-stage formulation by introducing a piecewise linear approximation of the value of the terminal inventory, to mitigate the two-stage myopic behavior. We compare the out-of-sample performance of the resulting models by rolling horizon simulations, within a data-driven setting, characterized by seasonality, bimodality, and correlations in the distribution of end item demand. Computational experiments suggest the potential benefit of adding a terminal value function and illustrate interesting patterns arising from demand correlations and the level of available capacity. The proposed approach can provide support to typical MRP/ERP systems, when a two-level approach is pursued, based on master production and final assembly scheduling.
\end{abstract}

\history{This is an Author’s Original Manuscript of an article published by Taylor \& Francis in the {\bf International Journal of Production Research} on 21.11.2023, available online: \url{https://doi.org/10.1080/00207543.2023.2283570} \\ Any citation must refer to its published version.}

\section{Introduction}
\label{sec:intro}

In the manufacturing domain, demand uncertainty is one of the main contributors to the difficulty of production planning problems, alongside production yield, quality, and machine reliability. A wide array of buffering tools has been devised to ease the difficulty of demand forecasting, including safety stocks \citep{Goncalves2020} and delaying assembly, in order to postpone product differentiation and take advantage of risk pooling effects (a well-known success story for this approach is the HP DeskJet case, see \cite{Lee1993}, \cite{Aviv2001}). The latter idea is 
exploited in 
{\em Assemble-To-Order} (ATO) manufacturing environments, where the long lead time to manufacture or procure components makes a pure Make-To-Order approach not feasible, so that components must be ordered under demand uncertainty. Yet, if final assembly is relatively fast, we do not need to stock end items, and final assembly can be carried out after demand for end items is observed. 

In the context of typical {\em Manufacturing Resource Planning}
(MRP) systems, as well as their descendants {\em Enterprise Resource Planning} (ERP) systems, this is reflected by a two-level approach.
A {\em Master Production Schedule} (MPS) is not defined for end items in the bills of materials, but rather for basic modules that are made to stock, according to forecasts of dependent demand, and then used according to a {\em Final Assembly Schedule} (FAS), driven by actual independent demand \citep{MasterSchedulingProud}.
However, traditional MRP/ERP systems may struggle with the complexity of finite capacity planning, which compounds with demand uncertainty. The typical approach is to introduce safety stocks, and to plan according to a rolling horizon strategy. The aim of this paper is to provide practical support to this kind of endeavor by prescriptive analytics, based on optimization models. We do so within a simplified setting, where bills of materials are flat and comprise two levels, basic modules/components, and end items, corresponding to the aforementioned MPS and FAS levels.

If a single assembly period is considered, the optimization problem can be formalized as a two-stage stochastic linear programming model, where at the first stage components (modules) are produced (under capacity constraints and subject to uncertainty on the demand of end items), while in the second one demand is observed and end items are assembled. Two-stage ATO models can be tackled with a limited computational effort \citep{Brandimarte2021} and may provide optimal solutions also in the multistage setting, if demand is \textit{independent and identically distributed}  (i.i.d.) across different time periods \citep{Gerchak1986}. Nevertheless, a longer planning horizon is necessary when the i.i.d.\ assumption does not apply, possibly due to predictable variability in demand, like a seasonal demand spike that cannot be met by the available production capacity and requires careful planning of inventory buffers. In this setting, the power and flexibility of multistage {\em Stochastic Programming} (SP) models may play a significant role, but their application is hindered by two issues: the lack of reliable information about probability distributions of the underlying risk factors and the exponential growth of the scenario tree. While the first issue can be tackled by adopting a data-driven approach, the second one could actually prevent the adoption of such techniques. 
Indeed, the literature shows that a considerable research effort has been devoted to proper scenario tree generation, in terms of scenario sampling and 
choice of its structure, i.e., breadth (branching factors) and depth (number of periods). In particular, limiting the number of periods may greatly reduce the size of the tree, at the price of a potential adverse effect in terms of myopia.

The nasty effects of truncated decision horizons are well-known even in a multiperiod deterministic setting. In \cite{Grinold1983}, issues related with the truncation of infinite-horizon problems are thoroughly analyzed. 
The addition of a terminal value is proposed by \cite{Fisher2001} to overcome end-of-horizon effects. A terminal state value is also a key ingredient in stochastic dynamic programming \citep{brandiDP}, and the integration of SP with stochastic control is proposed by \cite{Konicz2015} for a financial application (see also \cite{ZiembaYasudakasai1998}). Nevertheless, while in the financial domain there is a remarkable body of knowledge concerning the value of terminal states for asset--liability management problems,
much less is available in the manufacturing domain. 

The contribution of this work is twofold.
First, we propose a computational methodology for approximating the value of residual components in inventory, adding an end-of-horizon term to the objective function, to counteract the myopia of models with short horizons. Second, we study the performance of different model formulations in an ATO context characterized by challenging demand features, such as bimodality, seasonality, and correlations, also considering the impact of available capacity for component production. We assume that only a limited amount of historical data is available.
More in detail, we compare the following families of models:
\begin{itemize}
    \item A plain two-stage stochastic linear programming model.
    \item A two-stage stochastic model enhanced with an estimate of the value of residual component inventory.
    \item A set of multi-stage 
    stochastic models characterized by different types of scenario trees, varying in breadth and depth.
    \item A set of deterministic models, buffering against demand uncertainty by different levels of safety stocks. 
\end{itemize}
We test these models in a realistic environment, by simulating the production system within a rolling horizon framework. 
Only the decisions in the first time period are implemented, and
the models are repeatedly solved after observing the new demand.   
Simulations are run out-of-sample, whereas the in-sample scenarios used in the models rely only on a fixed limited amount of demand observations.
Computational experiments show that including an estimate of the value of the residual inventory can greatly improve the performance of the two-stage model, making it the best among those we tested. Moreover, we show that considering a longer planning horizon can be detrimental in terms of performance, and not only from a computational time point of view.
 
Since the problem is characterized by several parameters, we present in this paper only the most interesting results that we have obtained. Nevertheless, the open-source code is available at \url{https://github.com/DanieleGioia/ATO}, to properly guarantee reproducibility and enable interested readers to carry out further experiments.

The paper is organized as follows. In Section \ref{sec:litRew} we review the literature about ATO problems and rolling horizon policies,
and position the paper within the context of approximate dynamic programming. 
In Section \ref{sec:modelForm} we describe the mathematical models used in the computational experiments. Then, in Section \ref{sec:generation} we introduce the general experimental setting, and in Section \ref{sec:results} we present the results of several computational experiments. Finally, Section \ref{sec:conclusions} discusses conclusions and future research directions.

\section{Literature review and paper positioning}\label{sec:litRew}

Some early references on ATO problems date back to the 80s \citep{Wemmerlov1984, Collier1982}. 
However, a relevant part of this literature is actually oriented 
towards continuous- or periodic-review inventory control models
within a supply chain management framework, 
where the aim is to replenish inventory and allocate components, 
rather than finite-capacity production planning. 
With the aim of building analytical models, a commonly adopted tool is stochastic modeling, often dealing with a single end item or just a very few. 
Examples of such approaches are \cite{Fang2008}, where a game-theoretic framework is used to investigate pricing and procurement strategies, and \cite{Fu2006}, where outsourcing is also considered
and an analytical framework is applied to characterize
properties of the optimal solution for a problem involving a single end item. 
Stochastic control is applied in \citep{PlambeckWard}, to cope with customer order sequencing.
We refer to \citep{Atan2017,Song2003} for a review of 
this stream of literature.

On the contrary, we tackle discrete-time, 
multi-item, finite-capacity ATO production planning problems.
The \emph{Bills of Materials} (BoM) are flat and include only two levels, end items and components (modules), which correspond to the two-level MPS/FAS approach, common in MRP systems. 
We assume components are manufactured under finite capacity constraints and final assembly is not a bottleneck, so we do not consider assembly capacity and cost. 
Demand is uncertain and backordering is not possible. Hence, unsatisfied demand is lost, with an adverse impact on the overall objective, which we assume is to maximize expected profit.

In more recent years, the literature on deterministic mathematical programming models for finite-capacity production planning includes papers dealing with demand uncertainty, where both stochastic programming
and robust optimization are applied. A recent review is \cite{MultistageLS}.
While some works consider single-level production planning, others deal with multi-level problems, where product structures are represented by bills of materials. 
However, ATO problems feature a different information structure, 
as we may make assembly decisions {\em after} observing demand.

Since in this paper we focus on stochastic programming models with recourse, 
we first discuss the related literature. Then, we provide some references 
on deterministic approaches based on safety stocks, which is more akin to
standard practice in MRP systems and is used as a benchmark in our 
computational experiments. Finally,
we present our work within the more general framework of sequential decisions under uncertainty and approximate dynamic programming.

\subsection{Approaches based on stochastic programming}\label{sec:SP_ATO}

Most of the literature concerning the application of stochastic programming models to finite-capacity production planning under demand uncertainty deals with extensions of classical lot-sizing models, possibly involving multiple levels. \cite{DolguiStochLotSizing} provide a review of non-deterministic lot sizing models.

One of the first papers applying stochastic programming on ATO system
is \cite{Jonsson1993}, where the problem is formulated as a stochastic integer programming model. Fairly small problem instances are solved
by progressive hedging and scenario aggregation. 
After almost 30 years of improvements in commercial 
solvers, even some realistic size instances can likely be solved to near optimality by off-the-shelves software, at least in a two-stage case. 
Two-stage models are tackled by \cite{Brandimarte2021} and then extended by \cite{Brandimarte2022} to introduce risk aversion. 
Per se, two-stage models
are limited to newsvendor-like models.
However, they may be applied within approximation and decomposition schemes. 
In fact, truly multistage models are quite difficult to deal with, due to the explosion of the scenario tree.
For instance, \cite{Englberger} 
deal with a multi-period problem, but they apply 
a two-stage approach by freezing some decisions made at the first stage.
In the literature, this approach leads to static-dynamic approaches, which have also been applied to production planning involving the assembly of components.
For instance, \cite{TheveninSDDP} apply a static-dynamic approach, with an emphasis on solution methods based on stochastic dual dynamic programming. 
\cite{TheveninMRP} adopt a similar modeling framework, but focus on the application of practical heuristics and
testing by simulation.
Specific solution approaches for ATO systems are adopted by \cite{PrimalDual}, who apply newsvendor-based decomposition and rounding approaches; they apply sophisticated solution and analysis tools, whereas the demand structure is rather simple. 
\cite{DolguiChanceC} apply chance-constrained stochastic programming, as they place emphasis on the replenishment of components from suppliers, under random lead times.
\cite{Jaarsveld2015}, too, deal with component replenishment and allocation, rather than finite-capacity production. See also \cite{Huang2015}.
\cite{Nonas2009} adopts stochastic programming, providing analytical insights for small-scale component replenishment problems, limited to just a few  (three) end items.

Unlike, e.g.,  \cite{PrimalDual} and  \cite{TheveninSDDP}, we deliberately steer away from ad-hoc solution approaches, as they may be fragile when dealing with general models.
Since we aim at providing support to two-level master production scheduling, we prefer to use standard solvers. Furthermore, we insist on complex demand patterns.
We allow for seasonality and deal with realistic issues related to ATO system, where some end items may be variations on a basic family. Hence, we shall cope with intra-family demand, which needs a better approach than just pairwise correlations. Differently from \cite{PrimalDual} and \cite{DolguiChanceC}, we do not consider component replenishment from suppliers, but finite capacity production. This results in challenging issues, especially when there are demand peaks, possibly due to seasonality, that must be suitably addressed.

Our problem setting is more related with, e.g., \cite{Englberger} and \cite{TheveninMRP}. However, we consider a different information structure, related to two-level, rather than single-level master production scheduling. On the other hand, we disregard lot sizing issues associated with setup costs and times.

We do not use sophisticated scenario generation strategies as in \cite{TheveninSDDP}, like quasi-Monte Carlo sampling, as we rely on a data-driven approach whereby probability distributions are unknown and only a few data are available. Furthermore, we do not evaluate performance in terms of in-sample computational efficiency, but rather in terms of actual expected profit estimated by realistic out-of-sample, rolling horizon simulations (where the true demand distributions are used).

\subsection{Approaches based on deterministic models 
and safety stock buffers}\label{sec:safetystocks}

The practical strategy adopted in typical MRP systems relies on deterministic planning, integrated with the use of safety stocks as a hedging tool, within a rolling horizon process. 
This is reflected in academic research based on the application
of simulation models and optimization heuristics to set 
safety stocks.

The problem of sizing safety stocks has been considered in several papers and a recent literature review is provided by \cite{Goncalves2020}, who 
analyze a set of 95 articles published from 1977 to 2019 and show that 65\% of them apply analytical techniques, often within an inventory control framework, and do not contain references to realistic applications. 
On the contrary, simulation seems to play a key role in practice, often with reference to MRP systems.
Simulation-based optimization is adopted by \cite{Gansterer2014} and 
\cite{SimheurSafety}, 
as well as by \cite{SafetyPPC} to evaluate the performance of MRP systems. 
A deterministic optimization model is proposed by 
\cite{uzsoy}, using safety stocks to buffer against demand uncertainty. They use simulation to assess actual performance 
within a rolling horizon framework, which is what we
also do in our computational experiments.

\subsection{Paper positioning and relationships with approximate dynamic programming}

The distinguishing feature of our work is the integration of stochastic programming models for ATO with tools to ease \emph{end-of-horizon effects} \citep{Grinold1983,Fisher2001}. 
We should note that this kind of issue has been addressed in the financial domain \cite{ZiembaYasudakasai1998,Konicz2015} too. However, in finance we essentially have to deal with one commodity, namely, wealth.
In an ATO environment, the matter is more complicated, due to the relationships among components used for the same end item. 
Finding an exact expression of the terminal value function is impractical. Nevertheless, experience with rollout algorithms \citep{BertsekasRollout} shows that even a rough approximation may be remarkably effective in improving performance. A useful strategy for this aim is to resort to separable approximations as proposed, e.g., by \cite{PowellSimao2009}, which is what we pursue in this paper.

The comparison of solutions based on different scenario tree structures has been tackled for financial portfolio choice problems by, e.g., \cite{BirgeTree, Blomvall2006}.
However, financial and manufacturing domains are deeply different. While in the financial domain care must be taken to avoid building a tree allowing arbitrage opportunities, which places additional requirements on the shape of the scenario tree, in the production field there are no such problems. Moreover, in finance, there is an abundance of data, which is a much scarcer commodity in manufacturing, 
requiring a data-driven approach.
These differences motivate a specific investigation within a manufacturing setting.

In principle, the problem that we address in the paper could be solved exactly
by stochastic dynamic programming, which is 
the reference approach to cope with sequential decisions under uncertainty. 
It is well-known, however, that 
the curse of dimensionality precludes its application to most
practical size problems. Nevertheless, concepts from 
dynamic programming may be fruitfully applied to devise high-quality 
approximate solution strategies. 
\cite{powell2021} classifies different key approaches relying, e.g.,
on value function or policy function approximations, cost function approximations, and lookahead strategies, which may be useful to put alternative approaches to ATO production planning into a common perspective. 
Standard approximate dynamic programming approaches relying purely on value or policy function approximations do not look promising. 
On the one hand, the interactions among different components through the bills of materials
result in a complicated value function. On the other one, a policy function 
mapping the system state into component production decisions should cope with capacity constraints coupling different components. 
Furthermore, pure multistage stochastic programming models, which introduce a lookahead by a scenario tree are limited by the exponential growth of the tree. 
Hence, the essential contribution of our paper is to integrate
a lookahead approach, based on scenario trees, with a value function approximation.  
As discussed by \cite{AlphaZeroBertsekas}, this is the key idea in many successful approaches to dynamic decisions. 
A final strategy categorized by \cite{powell2021} is the approximation of a stochastic problem by a deterministic one, introducing buffers to hedge against uncertainty. In order to devise a competing approach, we do so by using safety stocks. 
%

\section{Model formulations} \label{sec:modelForm}

In this section, we introduce our model formulations of  
the ATO problem. Let $\mathcal{I} = \{1,\dots,I\}$ be the set of components, $\mathcal{J} = \{1,\dots,J\}$ the set of end items, and $\mathcal{M} = \{1,\dots,M\}$ the set of production resources (e.g., machines). Moreover, let us define the following parameters:
\begin{itemize}
    \item $C_i$ cost of component $i\in\mathcal{I}$.
    \item $P_j$ selling price of end item $j \in\mathcal{J}$.
    \item $K_j$ lost sale penalty of end item $j \in\mathcal{J}$.
    \item $H_i$ inventory holding cost of component $i\in\mathcal{I}$.
    \item $L_m$ production availability (time) of machine $m \in\mathcal{M}$.
    \item $T_{im}$ time required to produce one component $i\in\mathcal{I}$ on machine $m\in\mathcal{M}$.
    \item $G_{ij}$ amount of component $i\in\mathcal{I}$ required for assembling  end item $j \in\mathcal{J}$, commonly known as \emph{gozinto factor} when describing a BoM.
    \item $\bar{I}_i^0$ initial inventory of component $i$.
\end{itemize}

In order to formulate a stochastic model, we have to represent demand uncertainty. 
The standard choice in the SP literature is a scenario tree. Following the notation in \cite{Brandimarte2006}, we define:
\begin{itemize}
    \item $\mathcal{N}$, the set of nodes of the scenario tree and $\mathcal{N}^+ = \mathcal{N} \setminus \{0\}$. Node $0$ is the root node, where here-and-now decisions are made.
    \item $p(n)$, the parent of node $n \in \mathcal{N}^+$. 
    \item $\pi^{[n]}$, the unconditional probability of node $n$ ($\pi^{[0]}=1$).
    \item $d_j^{[n]}$, the demand for end item $j$ at node $n \in \mathcal{N}$.
\end{itemize}
We define as \emph{branching factor} the number of children (immediate successor nodes)
of each node at a given level of the tree. For example, by repeating three times a branching factor 2, the vector $[2,2,2]$ identifies a binary tree over four time periods, including the current time corresponding to the root node, which results in eight scenarios. A scenario is a sequence of nodes from the root node to a leaf node, i.e., a realization of the underlying stochastic process. All of the nodes at the same level are associated with
the same time instant.
For each node $n \in \mathcal{N}$, we define the following decision variables (restricted to non-negative integers):
\begin{itemize}
\item $x_i^{[n]}  \in \mathbb{Z}^+$, the amount of components $i\in\mathcal{I}$ produced at node $n \in \mathcal{N}$.
\item $I_i^{[n]} \in \mathbb{Z}^+$, the available inventory of components $i\in\mathcal{I}$ at node $n \in \mathcal{N}$.
\item $y_j^{[n]} \in \mathbb{Z}^+$, the amount of end items $j\in\mathcal{J}$ assembled at node $n\in \mathcal{N}$.
\item $l_j^{[n]} \in \mathbb{Z}^+$, the lost sales of end item $j\in\mathcal{J}$ at node $n\in \mathcal{N}$.
\end{itemize}

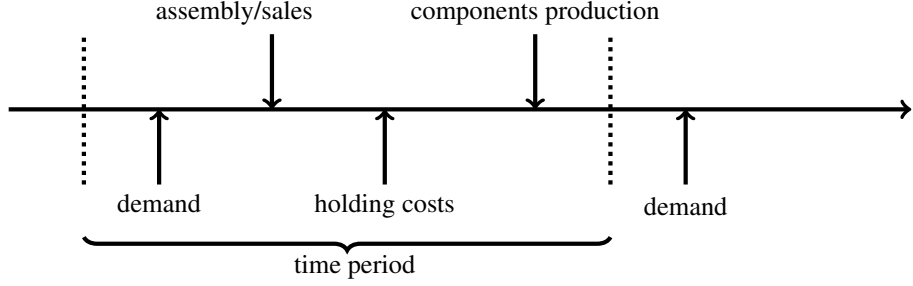
\begin{figure*}[tb]
    \centering
    \begin{tikzpicture}
    \draw[ultra thick, ->] (-2,0) -- (10,0);

    \draw[ultra thick,dotted] (-1,-1) -- (-1,1);
    
    \draw[ultra thick] (0,-1) node[below=0pt,thick] {demand} node[above=4pt] {};
    \draw[ultra thick, ->] (0,-1) -- (0,0);
    
    \draw[ultra thick] (1,1) node[above=0pt,thick] {assembly/sales} node[above=4pt] {};
    \draw[ultra thick, ->] (1.5,1) -- (1.5,0);
    
    \draw[ultra thick] (3,-1) node[below=0pt,thick] {holding costs} node[above=4pt] {};
    \draw[ultra thick, ->] (3,-1) -- (3,0);
    
    \draw[ultra thick] (5,1) node[above=0pt,thick] {components production} node[above=4pt] {};
    \draw[ultra thick, ->] (5,1) -- (5,0);

    \draw[ultra thick,dotted] (6,-1) -- (6,1);

    \draw[ultra thick,decorate, decoration={brace, amplitude=5pt}] (6,-1.7) -- (-1,-1.7);
    \draw[ultra thick] (2.6,-1.7) node[below=3pt,thick] {time period} node[above=4pt] {};
    
    \draw[ultra thick] (7,-1) node[below=1pt,thick] {demand} node[above=4pt] {};
    \draw[ultra thick, ->] (7,-1) -- (7,0);

    \end{tikzpicture}
    \caption{Event sequence at each node in the scenario tree.}\label{fig:informationFlow}
\end{figure*}

It is important to clarify the event sequence taking place at each discrete time instant. 
At each node, we first observe the new demand for end items. Then, we satisfy the demand using the components available in the inventory,
and we pay holding costs for components 
and (possibly) lost sales costs for end items. Finally, new components are produced. The flow of events is represented in Figure~\ref{fig:informationFlow}.
Understanding these dynamics will play a major role in the rolling horizon
simulations, and this is why we also consider end-item demand at the root node.  

The resulting model is:

\begin{alignat}{3}
&\!\!\max_{x,y} \quad &&  \sum_{n \in \mathcal{N}} \pi^{[n]} \!\left[ \, \Js (P_j y_j^{[n]} - K_j l_j^{[n]}) - \Is ( C_i x_i^{[n]} \right.&& \left. + H_i I_i^{[n]} )\right]  \label{eq:of_multi}\\
&\st \quad && \Is T_{im} x_i^{[n]} \le L_m  &&\ForM, \ForN \label{eq:prodConst_multi} \\
& && I_i^{[n]} = I_i^{[p(n)]} +  x_i^{[p(n)]} 
- \Js G_{ij} y_j^{[n]} && \ForI, \ForN^+ \label{eq:inventoryConstr_multi} \\
& && I_i^{[0]} = \bar{I}_i^0 - \Js G_{ij} y_j^{[0]} && \ForI \label{eq:initialInventory_multi} \\
& && y_j^{[n]} + l_j^{[n]} = d_j^{[n]} \label{eq:lostSales_multi} && \ForJ, \ForN \\
& && x_i^{[n]}, I_i^{[n]} \in \mathbb{Z}^+ &&\ForI, \ForN \label{eq:nonNegXI_multi} \\
& && y_j^{[n]}, l_j^{[n]} \in \mathbb{Z}^+ &&\ForJ, \ForN. \label{eq:nonNegLY_multi}
\end{alignat}
The objective function \eqref{eq:of_multi} is the expected net profit, expressed as the expected revenue, minus lost sales penalties and inventory holding costs, minus the cost of components production. Constraints \eqref{eq:prodConst_multi} represent the capacity constraints on component production. Constraints \eqref{eq:inventoryConstr_multi} describe inventory balance of components and define the ATO dynamics. In particular, the inventory at node $n$ is equal to the inventory at the parent node $p(n)$, plus the produced amount
at node $p(n)$, minus the quantity used for assembly at node $n$. Constraints \eqref{eq:inventoryConstr_multi}, together with non-negativity of inventory, preclude assembling end items for which there are not enough components. Note that, consistently with to the
dynamics that we have described before, 
only items that were available in inventory or  
were produced at the parent node $p(n)$ can be used for assembly at node $n$.
Initial conditions at the root node are imposed by constraints \eqref{eq:initialInventory_multi}. Constraints \eqref{eq:lostSales_multi} link end item assembly and lost sales with demand at each node $n$. Finally, constraints \eqref{eq:nonNegXI_multi} and \eqref{eq:nonNegLY_multi} specify that the decision variables
are non-negative integers.

\begin{figure}
    \centering
    \includegraphics[width=0.85\textwidth]{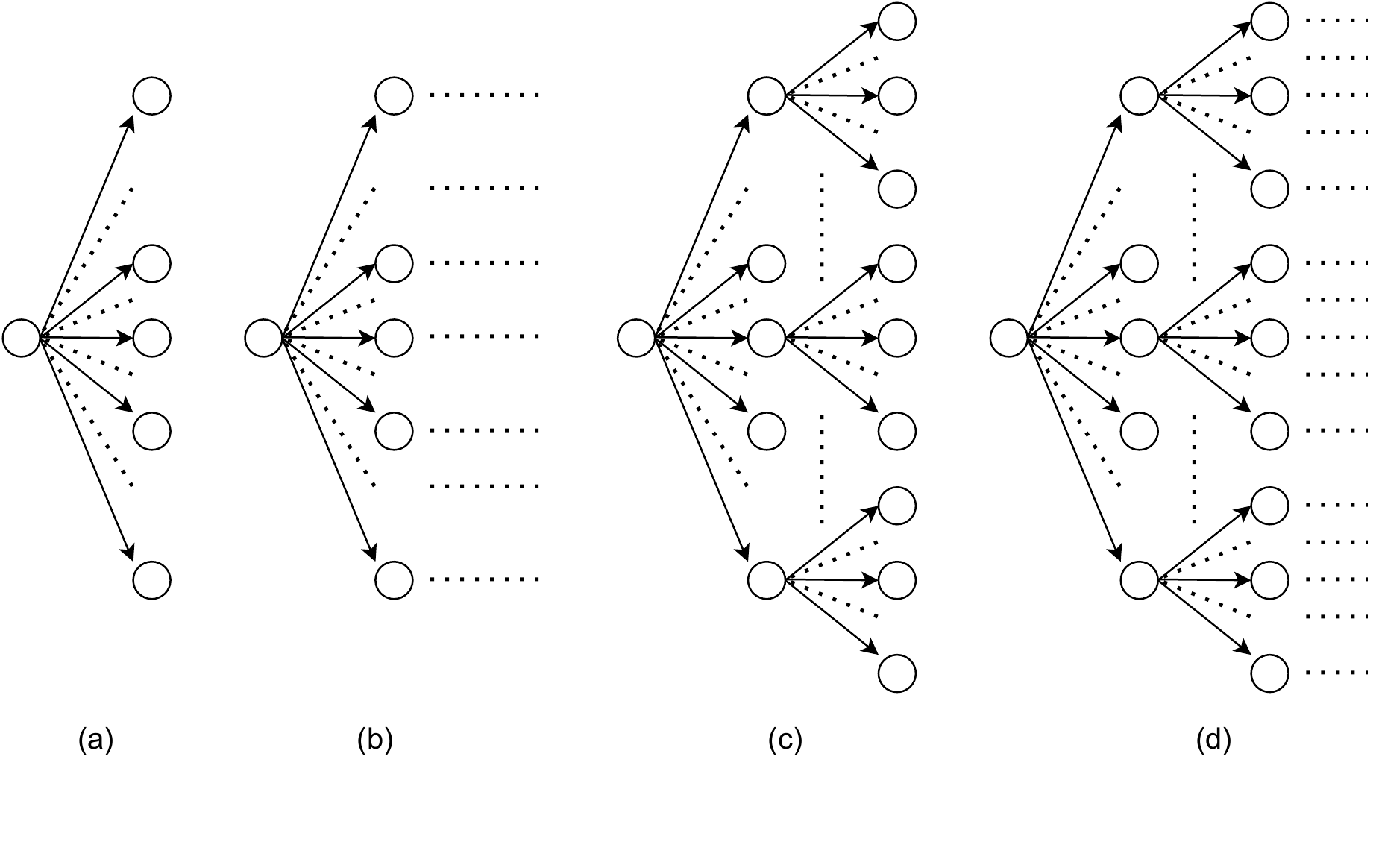}
    \caption{ Different scenario trees adopted: (a) \TS and \FOSVA, (b) \MP, (c) \MSt, (d) \MStp.}
    \label{fig:stages}
\end{figure}

The size of the model \eqref{eq:of_multi} - \eqref{eq:nonNegLY_multi}, increases exponentially as more time periods and greater branching factors are considered.
In the following, we explore reformulation strategies to solve it when a limited amount of data is available. In particular, we consider the following strategies:
\begin{itemize}
    \item a plain two-stage formulation;
    \item a two-stage formulation that associates a value to the terminal inventory state;
    \item alternative multistage formulations characterized by different branching factors;
    \item a deterministic model, where uncertainty is hedged by safety stocks.
\end{itemize}
We start with the shallowest tree, consisting of just two stages as shown in Figure \ref{fig:stages}(a). In this setting, for notational convenience, we drop the superscript $*^{[0]}$ for variables associated with the root node, and use the notation $*^s$ instead of $*^{[n]}$ for second-stage variables. Moreover, we rename $\mathcal{N}^+$ as $\mathcal{S}$, consistently with 
the two-stage programming literature. With these changes, model \eqref{eq:of_multi} - \eqref{eq:nonNegLY_multi} becomes:

\begin{alignat}{3}
&\max_{x,y}\quad && \Js (P_j y_j - K_j l_j) - \Is (H_i I_i + C_i x_i) + \Ss \pi^s && \left[ \Js (P_j y_j^s - K_j l_j^s)
   - \Is H_i I_i^s \right] \label{eq:of_ts}\\
&\,\st\quad && \Is T_{im} x_i \le L_m  &&\ForM \label{eq:prodConst_ts}\\
&& & y_j^s + l_j^s = d_j^s  &&\forall j \in \mathcal{J} , \forall s \in \mathcal{S} \label{eq:lostSales_ts}\\
&& & y_j + l_j = d_j  &&\forall j \in \mathcal{J} \label{eq:lostSales0_ts}\\
&& & I_i^s = I_i +x_i  - \Js G_{ij}y_j^s   &&\ForI ,\ForS \label{eq:inventoryConstr_ts}\\
&& & I_i = \bar{I}_i^0 - \Js G_{ij}y_j   &&\ForI \label{eq:inventoryConstr0_ts}\\
&& & x_i, I_i, I_i^s \in \mathbb{Z}^+  &&\ForI \label{eq:nonNegXI_ts} \\
&& & l_j^s, y_j^s, l_j, y_j \in \mathbb{Z}^+   &&\ForJ ,  \ForS
\label{eq:nonNegLY_ts}
\end{alignat}

In the following, we refer to model \eqref{eq:of_ts} - \eqref{eq:nonNegLY_ts} as \emph{two-stage} (\TS). \TS can be solved by off-the-shelves software with a 
fairly large number of scenarios \citep{Brandimarte2021}. Nevertheless, the optimal solutions to this problem are myopic. Specifically, for each component produced but not assembled, the inventory cost $H_i$ is paid and no gain is received. Thus, 
there are no explicit incentives to leave components in inventory after assembly; leftover inventory will just be the result of hedging demand uncertainty in the immediate future (the second stage). 

In order to overcome the myopic behavior of \TS, we consider a model that associates a value to the terminal inventory state. In other words, we add an extra term to the objective function to approximate the value of the components that are produced but are not immediately used for assembly in the second stage. By doing so, we obtain the following second-stage contribution:
\begin{alignat}{3}
\Ss \pi^s 
\!\left[\, \Js (P_j y_j^s - K_j l_j^s)\right. && \left .- \Is H_i I_i^s + V(I_1^s, \dots, I_I^s) \right] \label{eq:of_fosva_raw}
\end{alignat}
subject to the same constraints 
\eqref{eq:prodConst_ts} - \eqref{eq:nonNegLY_ts} of the \TS model formulation,
where $V(I^s)=V(I_1^s, \dots, I_I^s)$ is a function measuring the value of the terminal state.

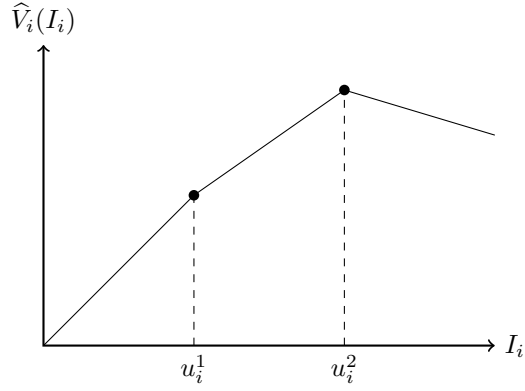
\begin{figure}[t]
    \centering
\begin{tikzpicture}[scale=2]
    \draw [<->,thick] (0,2) node (yaxis) [above] {$\widehat{V}_i(I_i)$}
        |- (3,0) node (xaxis) [right] {$I_i$};
    \draw[black] (0,0) coordinate (a_1) -- (1,1) coordinate (a_2);
    \draw[black] (1,1) coordinate (b_1) -- (2,1.7) coordinate (b_2);
    \draw[black] (2,1.7)  coordinate (c_1) -- (3,1.4) coordinate (c_2);
    \coordinate (c) at (intersection of a_1--a_2 and b_1--b_2);
    \coordinate (d) at (intersection of b_1--b_2 and c_1--c_2);
    \draw[dashed] (1,0) node[below] {$u^1_i$} -- (1,1) ;
    \draw[dashed] (2,0) node[below] {$u^2_i$} -- (2,1.7) ;
    \fill[black] (c) circle (1pt);
    \fill[black] (d) circle (1pt);
\end{tikzpicture}
\caption{ Example of a concave approximate value function.}
\label{fig:valueExample}
\end{figure}

Theoretically, one way of obtaining this value function is by 
infinite-horizon dynamic programming \citep{schrage}.
Since finding an exact and tractable expression of value functions is generally not feasible \citep{brandiDP}, we devise
a suitable approximation based on the following ideas:
\begin{enumerate}
    \item The value function can be approximated by
    decomposition with respect to components, 
    i.e., by setting
    $V(I) \approx \sum_{i\in \mathcal{I}} \widehat{V}_i(I_i)$.
    Clearly, by doing so we neglect interactions among components. Nevertheless, this decomposition approach is quite common and works often fairly well in practice  \citep{Powell2011}.
    \item To preserve the linearity of the model, each term of the sum is approximated by a piecewise linear function. 
    \item To avoid an excessive computational burden, the information needed to build the piecewise linear function is collected by solving the two-stage model \TS for different values of the initial inventory state. In other words, we assume $V(I) \approx \TS(I)$, where $\TS(I)$ denotes the optimal value of the two-stage model obtained by setting  
    $\bar{I}^0 = I$ in Eq \eqref{eq:inventoryConstr0_ts}.
\end{enumerate}
Putting all of these ideas together, we approximate the value of terminal component inventories as a sum of piecewise linear concave functions like the one depicted in  Figure~\ref{fig:valueExample}.
Concavity of the value function is a consequence of the fact that it is the optimal value of an LP in maximization form, as a function of the right-hand sides of constraints 
\eqref{eq:inventoryConstr0_ts} [see, e.g., \cite[Theorem 4.38]{GobernaBookNonlinear}]. Each individual term in the approximating sum is concave, as concavity is preserved by considering a concave function along its coordinates \citep{BoydBookConvex}.
As shown in Figure~\ref{fig:valueExample}, a slope can be negative for excessive inventory levels. It is important to observe that available inventory does not contribute to production costs in the \TS model, but there is an inventory leftover cost, which
penalizes huge inventories. These are further discouraged by the inventory holding cost in Eq \eqref{eq:of_fosva_raw}.
The term $\widehat{V}_i(I_i^s)$ can also be interpreted as a markdown value of the components, which is included in a term like  $\eta_i(I_i^s) = \widehat{V}_i(I_i^s) - H_i I_i^s$. If $\eta_i(I_i^s)$ is positive, leftover components of type $i$ contribute to value, even if they are not sold within an assembled end item. This will be the case if future stages require a quantity of components that is beyond production capacity.

To define $\widehat{V}_i(\cdot)$, we need to select a sequence of $K$ breakpoints (denoted by $u^k_i$) and a sequence of corresponding slopes (denoted by $v^k_i$) such that: \begin{equation}\label{eq:piecewiseV}
    \widehat{V}_i(I_i) = \begin{cases}
    v^0_i I_i, \quad & 0\le I_i \le u^1_i\\
    v^1_i (I_i - u^1_i) + v^0_i u^1_i, \quad & u^1_i\le I_i \le u^2_i\\
    \vdots\\
    v^k_i (I_i - u^{k-1}_i) + \sum_{l=1}^{k-1} v^l_i ( u^{l}_i - u^{l-1}_i), \quad & u^{k-1}_i\le I_i \le u^k_i\\
    \vdots
    \end{cases}
\end{equation}
where we fix $u^0_i = 0$, $\ForI$, and the sequence of slopes is non-increasing due to concavity. To choose the set of breakpoints and slopes, we consider a technique similar to \emph{Concave Adaptive Value Estimation} (CAVE) introduced by \cite{Godfrey2001, Godfrey2002}. In order to avoid incurring the curse of dimensionality, breakpoints are not defined by a grid, but they 
are obtained by plain Monte Carlo sampling. At each iteration, we generate a point $\hat{I} = [\hat{I}_1, \dots, \hat{I}_I]$ by sampling each component inventory level from a uniform distribution between $0$ and an upper bound $I^{\max}$. Since we assume no inventory limit, $I^{\max}$ must be set empirically. A good rule of thumb is to set $I^{\max}$ to a multiple of the average demand; this is a value large enough to ensure that, up to a certain degree of confidence, demand in the next stages is satisfied. 

Given a sampled breakpoint value $\hat{I}$, to build the piecewise linear approximation of Figure~\ref{fig:valueExample}, we have to generate information about the slope on its right and on its left. Thus, we apply a straightforward scheme based on \emph{finite
differences}. The left and right slopes of $\widehat{V}_i$ 
at point $\hat{I}_i$, denoted by $\pi^+_i$ and $\pi^-_i$, respectively,
are computed as:
\begin{equation}\label{eq:defPi}
    \pi^+_i = \frac{\TS(\hat{I} + \varepsilon e_i)-\TS(\hat{I})}{\varepsilon}, \qquad
    \pi^-_i = \frac{\TS(\hat{I}) -\TS(\hat{I} - \varepsilon e_i)}{\varepsilon},
\end{equation}
where $e_i$ is the $i$-th element of the orthonormal basis 
of $\mathbb{R}^I$, and $\varepsilon$ is 
the perturbation of the approximation step. 
The idea is to perturb the inventory of component $i$ and to check
the impact on the optimal value of model \TS, in terms
of a forward and a backward finite difference. 
We assume that the number of scenarios in the two-stage model is sufficient to obtain suitably reliable information. Nevertheless, some degree of smoothing and dynamic adaptation is needed when adding new sampled points. To this aim, we apply a strategy whose pseudo-code is given in Algorithm~\ref{alg:iCAVE}.

\begin{algorithm}[tb]
\footnotesize
\caption{First order stock value approximation}\label{alg:iCAVE}
\begin{algorithmic}
\State $u_i = [0], v_i=[0] \  \forall\ i$
\For{\texttt{count=1:K}}
    \State \texttt{generate a random point} $\hat{I} = [\hat{I}_1, \dots, \hat{I}_I]$ \texttt{with} $\hat{I}_i \sim \mathcal{U}[0, I^{\max}] \ \forall\ i$
    \For{\texttt{i=1:I}}
        \State \texttt{add} $\hat{I}_i$ \texttt{in the sorted array} $u_i$, \texttt{call the position pos\_i} 
        \State \texttt{update the vector} $v_i$ \texttt{adding in pos\_i the value} $v_i^{\texttt{pos\_i - 1}}$ 
        \State \texttt{estimate the left and right derivatives of $V(\cdot)$ in $\hat{I}_i$}: $\pi^+_i$, and $\pi^-_i$ 
        \State $\nu^+_i \gets (1-\alpha) v_i^{\texttt{pos\_i}} + \alpha \pi^+_i$
        \State $\nu^-_i \gets (1-\alpha) v_i^{\texttt{pos\_i}} + \alpha \pi^-_i$
        \For{\texttt{k=1:$\min[$pos\_i,len($v_i$)]}}
            \If{$v^k_i \leq \nu^+_i$} 
                \State $v^k_i \gets \nu^+_i $
            \EndIf
        \EndFor
        \For{\texttt{k=pos\_i:len($v_i$)}}
            \If{$v^k_i \geq \nu^-_k$} 
                \State $v^k_i \gets \nu^-_k $
            \EndIf
        \EndFor
    \EndFor
\EndFor
\end{algorithmic}
\end{algorithm}

The algorithm updates the vectors of slopes and breakpoints independently for each component. It inserts $\hat{I}_i$ in the vector of breakpoints ($u_i$) in ascending order, where \texttt{pos\_i} is the position of $\hat{I}_i$ in $u_i$. Then, it updates the vector of slopes ($v_i$) by inserting a new element in \texttt{pos\_i} equal to the element in position \texttt{pos\_i-1}, which yields vectors $u$ and $v$ with the same length. Moreover, since $v$ in each iteration must have non-increasing components, setting the element in position \texttt{pos\_i-1} also in position \texttt{pos\_i} makes sure that the sequence of values in $v$ is non-increasing. 
The smoothed values $\nu^+ = (1-\alpha) v_i^{\texttt{pos\_i}} + \alpha \pi^+_i$ and $\nu^-= (1-\alpha) v_i^{\texttt{pos\_i}} + \alpha \pi^-_i$ are used to update the vector $v$. Specifically, all of the values on the left of \texttt{pos\_i} smaller than $\nu^+$ are set to $\nu^+$, and all of the values on the right of \texttt{pos\_i} greater than $\nu^-$ are set to $\nu^-$. These two steps create a non-increasing sequence of values for $v$.
The algorithm refines the approximation $K$ times in order to generate a reasonable approximation. The choice of $K$ needs some tuning because, on the one hand,
the more points are added the better the approximation. On the other hand, however, new breakpoints generate new variables in the resulting optimization problem. 

Based on the proposed approximation, and using $\mathcal{K}$ to denote the set of breakpoints of the piecewise linear approximation $\hat{V}$, we obtain the following model:
\begin{alignat}{3}
&\max_{x,y}\quad &&\Js(P_j y_j - K_j l_j) - \Is (H_i I_i + C_i x_i) + \nonumber \\&  && \quad\quad\quad +\Ss \pi^s \left[\, \Js (P_j y_j^s - K_j l_j^s) - \right.&&\left. \Is (H_i I_i^s + \Ks v^k_j m^s_{ik}) \right] \label{eq:of_fosva}\\
&\,\st\quad && I_i^s = \Ks m_{ik}^s   && \ForI \,\,\ForS \label{eq:sumM_fosva}\\
&& &  m_{ik}^s \le  u^{k}_i - u^{k-1}_i  && \ForI \,\,\ForS \,\, \ForK \label{eq:pieceM_fosva}\\
&& & m_{ik}^s \ge 0   &&\ForI \,\, \ForS \,\, \ForK \\
&& & \eqref{eq:prodConst_ts}, \eqref{eq:lostSales_ts},
\eqref{eq:lostSales0_ts},\eqref{eq:inventoryConstr_ts},
\eqref{eq:inventoryConstr0_ts},\eqref{eq:nonNegXI_ts},\eqref{eq:nonNegLY_ts},\nonumber
\end{alignat}
where the auxiliary variables $m_{ik}^s$ represent the piecewise linear approximation of Eq.\ \eqref{eq:piecewiseV}.
We refer to this model formulation as \emph{First Order Stock Value Approximation} (\FOSVA). Both \FOSVA and \TS rely on the same type of scenario tree in Figure \ref{fig:stages}(a), but \FOSVA is expected to be less myopic, at the cost of a preliminary learning step.
This step requires an offline estimation of the coefficients by solving $3\times K\times I$ different \TS problems.

Another possibility to make the \TS model less myopic is to consider a scenario tree as the one in Figure~\ref{fig:stages}(b), which is two-stage and multiperiod. Specifically, we add nodes assuming a known demand after the first branching. If independence holds, it can make sense to generate such values considering the expected value of demand (or its estimate) for each node. The information structure of this model, which we will refer to as \MP (where \texttt{<n>} is the length of the terminal linear graph), aims at mitigating end-of-horizon effects by introducing a sort of certainty equivalent steady state, without increasing the number of scenarios. It is worth noting that this structure, which has proven to be beneficial in manufacturing applications \citep{Brandimarte2006}, cannot be used in financial applications since it would lead to arbitrage opportunities \citep{BirgeTree}.  

We also consider the scenario tree structures of Figures~\ref{fig:stages}(c) and \ref{fig:stages}(d). In particular, the scenario tree of Figure~\ref{fig:stages}(c),
labeled as \MSt, relies on a more detailed 
description of uncertainty in the near future, as it comprises three stages.
The structure of Figure~\ref{fig:stages}(d) adds to this structure a sequence of nodes with a branching factor equal to 1. This is an extension of \MSt along the lines of \MP, and it is labeled as \MStp (where \texttt{<n>} is the length of the terminal linear graph).

\begin{figure}
    \centering
    \includegraphics[width=0.5\textwidth]{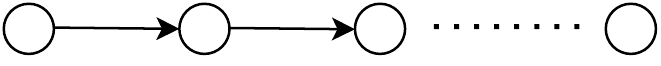}
    \caption{Scenario tree with estimates of the expected value of demand for each stage.}
    \label{fig:det}
\end{figure}

The last model formulation that we consider
aims at minimizing the computational effort, and relies on a sequence
of branching factors like 
$[1,1,\dots,1]$, as 
illustrated in Figure~\ref{fig:det}.
This is a deterministic model, where demand uncertainty is disregarded, and an estimate of expected demand is associated with each node.
The computational saving allows to consider a longer horizon, such as an entire year. 
However, ignoring demand uncertainty may lead to poor results, especially with complex distributions like the ones that we will consider in the computational experiments. 
A widely adopted strategy to hedge against uncertainty is to include a safety stock constraint into the model, which can be represented as follows:
\begin{equation}\label{eq:SSconstr}
     I_i^{[n]} \ge \sum_{j \in \mathcal{J}}G_{ij}d_j^{\,\alpha} \qquad \ForI, \ForN,
\end{equation}
where $d_j^{\,\alpha}$ represents the empirical $\alpha-$quantile of the demand distribution for end item $j$. Thus, the greater $\alpha$, the higher the safety stock. We refer to this family of heuristics as \SSq{$\alpha$}.

\section{Problem instance generation and design of experiments}\label{sec:generation}

A problem instance is characterized by a set of technological and economic data, by the probability distribution of demand, and by the structure of the scenario tree.
The generation of problem instances aims at supporting the computational experiments, where the proposed model formulations are compared on a set of out-of-sample rolling horizon simulations.
We create demand settings that are quite difficult to handle, in order to test the efficiency in using the available information by different models.  

\subsection{Technological and economical data}\label{sec:technoecodata}

From the technological and economic viewpoint, a problem instance is characterized by the following data
(see Table~\ref{tab:instanceParams}): 
\begin{itemize}
    \item capacity requirements and availability ($T_{im}$, $L_m$);
    \item economic parameters, i.e., costs, penalties, and sales prices  ($C_i$, $H_i$, $K_j$, and $P_j$);
    \item product structures, i.e., the gozinto factors ($G_{ij}$).
\end{itemize}

First of all, we set the number of components $I$, machines $M$, and end items $J$ according to Table~\ref{tab:instanceParams}. Afterward, we define a \emph{tightness factor} $\gamma$ to generate the available capacity for each machine $L_m$, 
with respect to the average capacity requirement. In formulas:
\begin{equation}
    L_m = \gamma \sum_{i\in\mathcal{I}} T_{im} \left(\sum_{j\in\mathcal{J}} G_{ij} \bar{d}_j\right),
\end{equation}
where $\bar{d}_j$ is the expected demand for end 
item $j$ that, given the complexity of the demand generation mechanism described below, is estimated by 
sampling $5000$ demand scenarios. In the experiments, we set $\gamma \in \{1.0, 1.1, 1.2, 1.3\}$, corresponding to instances where production capacity is equal to the average demand over the entire horizon ($\gamma = 1.0$), and up to 30\% more. The higher $\gamma$, the easier is to handle seasonal peaks, since it is possible to produce more components than those needed to satisfy average demand.

\begin{table}[t]
\small
\caption{Parameters defining a problem instance.}\label{tab:instanceParams}
\medskip
\begin{center}
\scalebox{0.95}{
\begin{tabular}{  l c  l c  }
Parameter & Value & Parameter & Value\\  
\hline
No.\ end items ($I$) & 35 & No.\ components ($J$) & 60\\
No.\ machines ($M$) & 5 & Tightness factor ($\gamma$) & \small{\{1, 1,1, 1.2, 1.3\}} \\
\hline
Low margin profit & \small{[0.05, 0.2]} & 
\% low margin items & 40\% \\
Medium margin profit & \small{[0.2, 0.4]} &
\% medium margin items & 30\% \\
High margin profit & \small{[0.4, 0.6]} 
& \% high margin items & 30\% \\
\hline
No.\ outcast items & 5 & No.\ families ($R$)& 5\\
No.\ items per family & \small{[12, 7, 5, 3, 3]} & No.\ 
common components & 2\\
No.\ components per family & \small{[11, 17, 12, 6, 9]} & &\\
\hline
Lost sale penalty $K_j$ & $0.2\cdot P_j$ & Holding cost $H_i$ & $0.1\cdot C_i $\\ 
\hline
\end{tabular}}
\end{center}
\end{table}
%

Concerning the economic values, we randomly generate the cost of components following a uniform distribution with range $[1,50]$. Then, we partition the end items into three different classes of profit margin: low, medium, and high. In Table~\ref{tab:instanceParams}, we define the proportion of items in each 
class and the ranges from which these margins are uniformly sampled. 
By considering the bill of material of each end item and the cost of required components,
we obtain the cost of each end item; the sales price is obtained by adding a profit markup, corresponding to each profit class.
By a similar token, we obtain inventory holding costs for components
and lost sales penalties for end items, as a percentage of costs and sales prices, respectively.

\begin{figure}[t]
\centering
\includegraphics[width=0.9\textwidth]{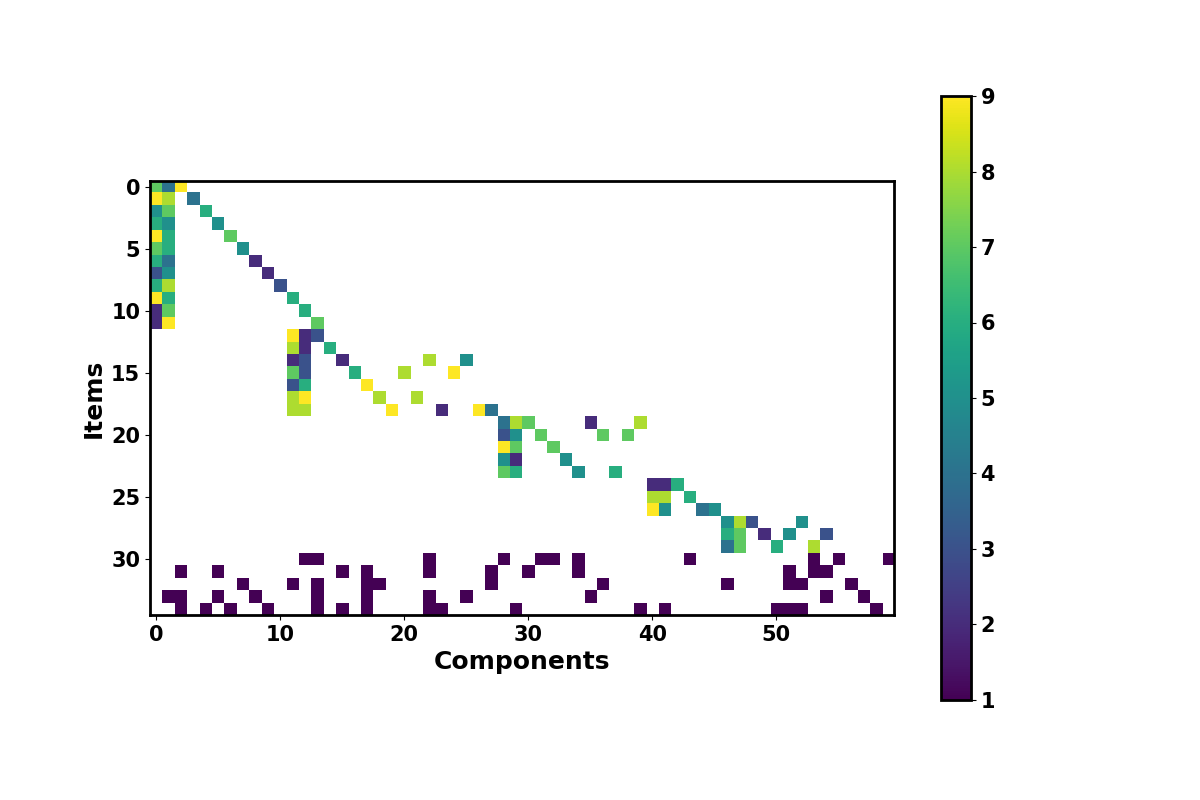}
\caption{A visualization of the gozinto matrix. The number of components required by each end item is represented by the color of each cell. The more yellowish the cells are, the higher the number of components.}\label{fig:gozinto}
\end{figure}

In the literature, several standard structures for the gozinto matrix have been proposed \citep{Atan2017}. Nevertheless, these structures are usually considered for their theoretical properties rather than their realism. We define the gozinto matrix as follows. First, we define families by partitioning the set of items: each family has a given number of end items and, within a family, end items have some required components, either specific or common with other members of the family. Moreover, we also introduce \emph{degenerate} families, consisting of a single end item. 
Such end items are, in a sense, apart from the others and called \emph{outcast items}. In Table~\ref{tab:instanceParams}, we report the characteristics of the \emph{standard} families (number of common vs.\ total components and number of items per family) and the number of outcast items. Once the total number of components of a (standard) family is assigned, we select the common components and the specific ones for each end item. Then, we sample how many pieces of each component are required by the end item, according to a discrete uniform distribution ranging from 1 to 9. On the contrary, outcast items do not share any specific pattern with others, and their components are picked according to a binomial distribution so that each component is included with probability $0.2$. 
The gozinto matrix, generated by the aforementioned procedure and adopted for the computational experiments, is depicted in Figure~\ref{fig:gozinto}. The diagonal blocks correspond to the standard families (five in this case), while the bottom rows correspond to the outcast items. Each standard family shares a number of common components (the first columns of each block).

\subsection{Demand distributions}\label{subsec:Demand}

Since we assume integer decision variables in our optimization models, we should adopt discrete probability distributions to generate demand scenarios.
Nevertheless, in order to take advantage of the richness of continuous distributions, we sample according to the continuous demand distributions that we describe in this section, and then we just round the obtained values. 

To emphasize the risk associated with production planning and the complications created by a distribution that is not easy to estimate empirically, especially in the case of a small number of observations, we consider seasonal patterns and bimodal distributions. Furthermore, we also investigate the impact of correlation patterns inside families. We consider monthly seasonality, modeled by factors that multiply both expected values and standard deviations of the base-case distributions. The multiplicative factors are given in Table~\ref{tab:seasonal}, and we apply them to all end items in the same way. The capacity tightness factors refer to the average (base-case) demand during the year. It is worth noting that with the maximum value of tightness ($\gamma = 1.3$), the maximum seasonal peak can be met, on average, without prior inventory.

%
\begin{table}[tb]
\caption{Demand pattern
for yearly multiplicative seasonality.}\label{tab:seasonal}
\medskip
\begin{center}
\begin{tabular}{ c c c c c c c c c c c c c}\hline
Month & 0 & 1 & 2 & 3 & 4 & 5 & 6 & 7 & 8 & 9 & 10 & 11 \\  
Seasonality & 1.0 & 1.1 & 0.9 & 0.8 & 1 & 0.8 & 1.2 & 1.3 & 1.2 & 1.0 & 0.8 & 0.9\\
\hline
\end{tabular}
\end{center}
\end{table}
%

Another relevant feature of end item demand under our assumptions is bimodality. 
Bimodal distributions may be observed in market conditions like:

\begin{itemize}
\item
\emph{New products}. It is often difficult to predict whether a newly launched product will become a top seller or a shelf warmer. 
This is typical of the fashion industry.
When a new collection is launched, it may be unclear which styles will prove to be popular. The corresponding demands are thus governed by multimodal distributions, where each mode reflects a particular popularity state.
\item
\emph{Big customers}. Suppliers serving many small customers with independent demands face an aggregate demand that is approximately Gaussian. In the presence of a large customer that accounts for a large share of the sales, however, demand may be multimodal due to irregular bulk orders.
\item
\emph{New market entrants}. The entrance (or exit) of a major competitor can have a significant impact on demand. 
\end{itemize}
We assume that distributions of all end items are bimodal for every time period. This is rather questionable in practice. In fact, if bimodality is due to uncertainty about the market penetration of a new item, after the first demand periods we should be able to figure out whether it is a success or a flop. However, our assumption aims at creating a difficult environment in which to compare model formulations. 

End items in different families can be completely unrelated in terms of consumer markets. Hence, we assume that there is no correlation between families. However, there can be a correlation within a family. If a consumer buys an end item with a certain configuration, she is not likely to buy a similar item with a different configuration, which implies a negative correlation. We sample from a bimodal distribution that arises by mixing, with 
mixing parameter $p$, two normal distributions with expected values $(\mu_1,\mu_2)$ and standard deviations $(\sigma_1,\sigma_2)$, respectively. Thus, when correlation is disregarded, we have
\begin{equation}\label{eq:BiModal}
    D^{\text{Uncorr}}_j \sim \mathsf{BN}(\mu_1,\mu_2,\sigma_1,\sigma_2,p), \qquad \ForJ,\footnote{We denote random variables with capital letters and their realizations with lowercase letters.}
\end{equation}
where $\mathsf{BN}$ refers to a bimodal arising by mixing two normal distributions.
We list the parameters of the employed bimodal distributions in Table~\ref{tab:deamandParams}.

%
\begin{table}[tb]
\caption{Bimodal Gaussian distribution parameters.}\label{tab:deamandParams}
\medskip
\begin{center}
\begin{tabular}{l c l c}
Parameter & Value & Parameter & Value\\  
\hline
$\mu_1$ & 300 & $\mu_2$ & 50 \\
$\sigma_1$ & 50 & $\sigma_2$ & 15\\
p & 0.8 & &\\
\hline
\end{tabular}
\end{center}
\end{table}
%

\begin{figure}[tb]
\centering
\includegraphics[width=0.7\textwidth]{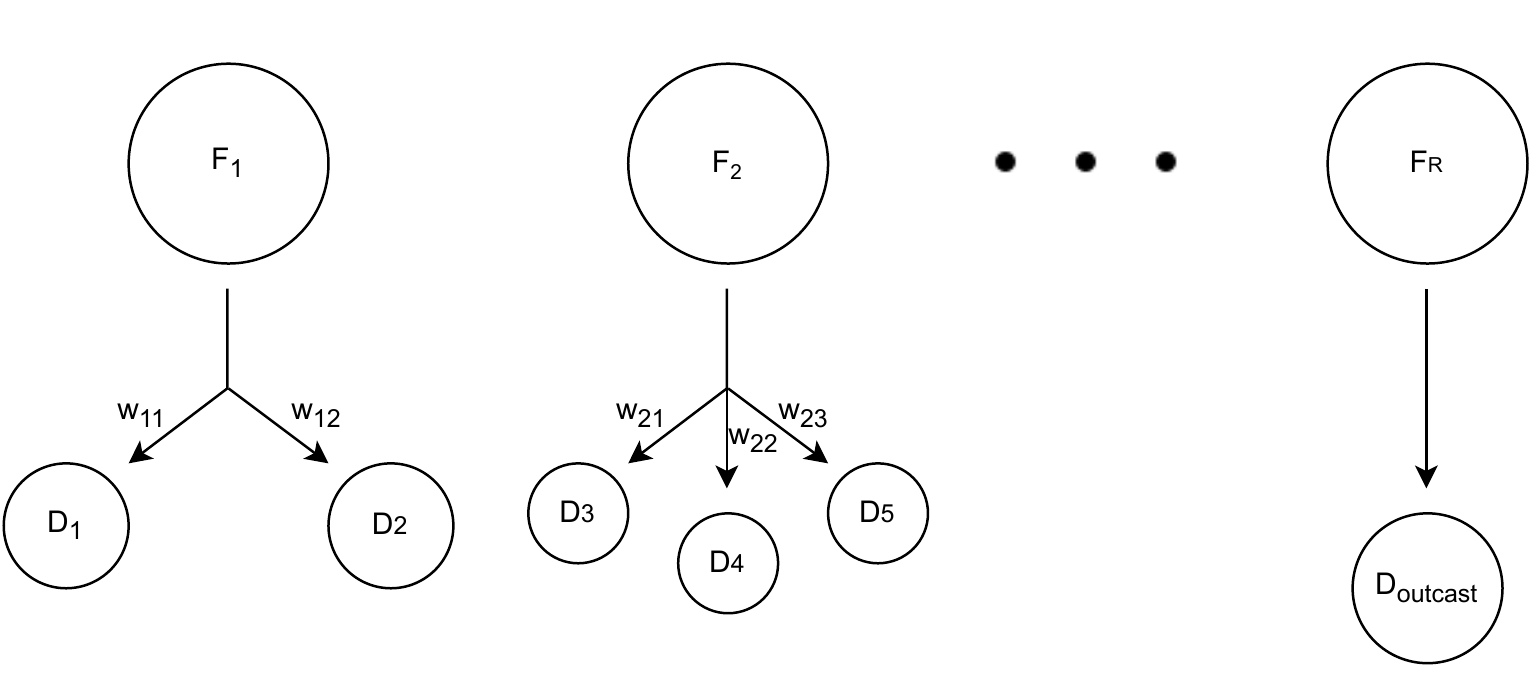}
\caption{Sampling demands with correlation within families.
Demands for families $F_1,F_2,\ldots,F_R$ are independently sampled.
Then, demands within families arise from disaggregation according to weights $w_{jr}$ (with the exception of 
outcast items). }\label{fig:dem}
\end{figure}

The case of correlation within families must be handled with some care, as arbitrary negative correlations may lead to a correlation matrix that is not positive semidefinite. 
Hence, we sample the correlated demand in two steps. We independently sample the aggregated demand $F_r$ for each family $r\in\{1,\dots, R\}$, according to a bimodal normal distribution, but with parameters $(n_r\mu_1,n_r\mu_2)$, $(\sqrt{n_r}\sigma_1,\sqrt{n_r}\sigma_2)$, and mixing parameter $p$, where $n_r$ is the number of end items in family $r$. Thus,
\begin{equation} \label{eq:biModalFam}
    F_r \sim \mathsf{BN}(n_r\mu_1,n_r\mu_2,\sqrt{n_r}\sigma_1,\sqrt{n_r}\sigma_2,p), \qquad \forall r \in\{1,\dots, R\}.
\end{equation}
Then, the overall demand of family $r$ is split among the items belonging to the family ($j \in J_r$), according to weights $w_{jr}$ randomly sampled from a \emph{Dirichlet distribution} with different parameters $\zeta_r$ for each family $r$, i.e.,
\begin{equation}
    W_{1r},\dots,W_{n_rr} \sim \mathsf{Dirch}(\zeta_r), \qquad \forall r \in\{1,\dots, R\}.
\end{equation}
This choice ensures that $\sum_{j \in J_r} w_{jr} = 1$, $w_{jr}>0,$ $\forall j \in J_r\, , r  \in\{1,\dots, R\} $. Clearly, in case of families consisting of just one item (as for the outcast items), the demand for the unique end item is equal to the demand of the family. The resulting distribution is
\begin{equation}
    D^{\text{Fcorr}}_j \sim w_{jr} F_r, \quad \forall j \in J_r, \quad \forall r \in\{1,\dots, R\}.
\end{equation}
In Figure~\ref{fig:dem}, we provide a graphical representation of the within-family correlation process. We report an example of the weights sampled from the Dirichlet distribution to generate the within-family correlated demands in Table~\ref{tab:expShareNominal}.

%
\begin{table}[tb]
\caption{Correlation weights per family.}\label{tab:expShareNominal}
\medskip
\begin{center}
\begin{tabular}{ c c } 
Family & Correlation weights  \\  
\hline
\multirow{2}{*}{1} & 8.45\%,  1.40\%,  5.63\%, 16.90\%,  5.63\%, 11.26\%, 14.08\%,  5.63\%,\\
& 8.45\%, 4.22\%,  7.04\%, 11.26\% \\
\hline
2 & 26.92\%,  3.84\%,  3.84\%, 19.23\%, 11.53\%,  7.69\%, 26.92\% \\
\hline
3 & 12.5\%, 25\%,  25\%,  12.5\%, 25\% \\
\hline
4 & 25\%, 25\%, 50\% \\
\hline
5 & 42.85\%, 14.28\%, 42.85\% \\
 \hline
\end{tabular}
\end{center}
\end{table}
%

\subsection{Scenario trees}

We consider time periods of one month and assume that 3, 5, or 10 years of previous demand data are available. We do not consider more than 10 years of available data, 
since demand may be subject to trends, items can become obsolete, etc. Thus, we have 36, 60, or 120 demand observations available for each end item, respectively. Nevertheless, given the monthly seasonality, we only have 3, 5, or 10 observations for each seasonal period. 
With such a small amount of data, estimating a demand distribution is out of the question, and one has to resort to purely data-driven approaches. 

Considering the entire seasonal pattern would require 12 time periods, which results in a scenario tree too big for commercial solvers. 
Moreover, considering horizons greater than 4 periods does not prove to be particularly useful, because the machine capacity does not allow for the production of components for such distant periods. We start with two-stage models (\TS, \FOSVA), in which only one future stage is considered. Then, we gradually increase the number of considered time periods to 2 (\MSt, \MPd), to 3 (\MStt, \MPt), and finally
to 4 (\MStq, \MPq). Following the branching strategy introduced in Section~\ref{sec:modelForm}, \MSt has two levels with all the available observations (i.e., we have branching factors of $[3,3]$, $[5,5]$, and $[10,10]$ for instances with 3, 5, and 10 years of data, respectively). \MStp and \MP have either one or two 
levels populated with all the available data. Both are followed by tails with nodes characterized by demands equal to the empirical mean (e.g., their branching factors are $[10,10,1,\dots,1]$ and $[10,1,\dots,1]$, respectively, assuming that 10 years of observations are available).

We also consider a two-stage model where scenarios are obtained by the entire set of demand observations, neglecting seasonality, and we denote this model as \TSn.
For example, in the case of 10 years of available observations, this model has 120 scenarios in its second stage. The role of \TSn is to provide a direct comparison with \FOSVA by proving that the value function estimation is not equivalent to considering all the observations together.


\section{Experimental results}\label{sec:results}


In this section, we compare the different scenario tree structures
based on rolling horizon, out-of-sample simulations. Both the in-sample data and the out-of-sample data used in the rolling horizon simulations are obtained from the same distribution. We consider a simulation horizon of two years, with monthly production decisions. The length of the test horizon allows us to appreciate the periodicity of the solutions and mitigate any dependence on the initial inventory, which is set to the real average demand for each component.

\subsection{Performance metrics}

The main metric, directly linked to the objective function of the optimization models, is  \emph{expected profit}. Profit, however, does not fully reveal certain characteristics of the solution behavior over the time horizon. 
Gaining further qualitative insights about how a specific model formulation behaves may provide additional and useful intuition in order to understand
the underlying mechanics and explain why an approach achieves good or bad performance.
Thus, we further investigate  \emph{average inventory level} (linked to holding costs) and  \emph{average lost sales}.
For expected profit and inventory level, we provide the percentage values with respect to the profit and the optimal inventory level of the problem with perfect information (Table~\ref{tab:profits},~\ref{tab:inventory}). However, comparing lost sales against the case of perfect information is not possible. In fact, since the perfect information model knows demand in advance, it has no lost sales on average. To overcome this problem, we provide the percentage values of the mean deviation of lost sales from the mean value of lost sales computed considering all of the methods, tightness configurations, and observed years in Table~\ref{tab:lostSales}. For example, a value of -10\% means that the method has ten percent fewer lost sales than the overall average quantity of lost sales. Conversely, a value of 100\% implies that the method has twice as many lost sales as the average value.
All the statistics are computed over 10 replications of the 24-month out-of-sample horizon.

Since solution times of the models generally do not exceed one minute, they are not reported and are out of the scope of the paper. 
It must be stressed that this is due to the limited size of the scenario tree, where we only use historical data according to a data-driven strategy. Indeed, the point of the paper is about using a limited amount of available data in the best way, not about solving large-scale stochastic programs. 
As we have discussed, in our study it is more relevant to check the actual solution quality by out-of-sample, rolling horizon simulations.

\subsection{Two-stage models}

\begin{figure}[t]
\centering
\begin{subfigure}[]{0.49\textwidth}{\includegraphics[width=\textwidth]{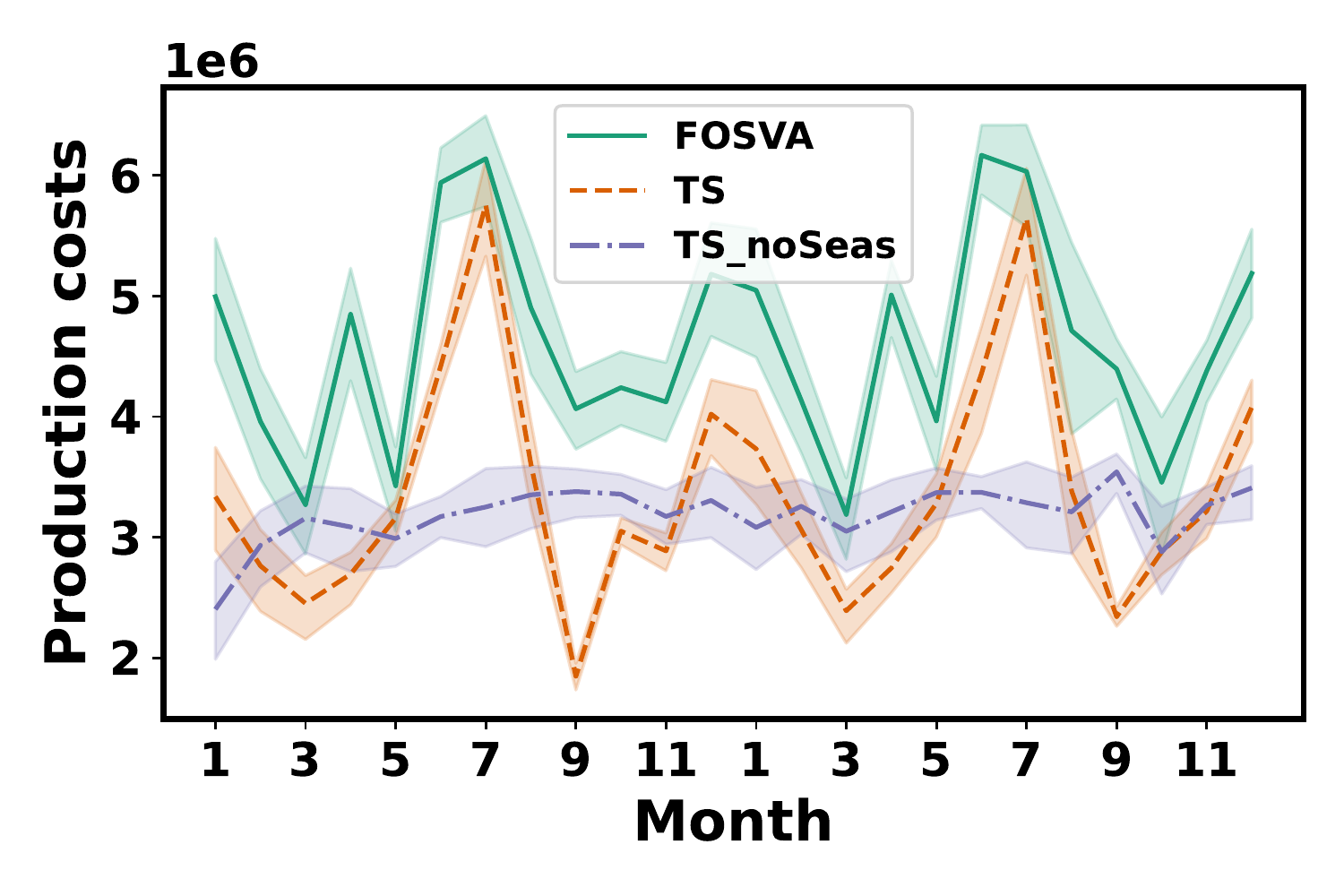}}
\end{subfigure}
\,
\begin{subfigure}[]{0.49\textwidth}
\includegraphics[width=\textwidth]{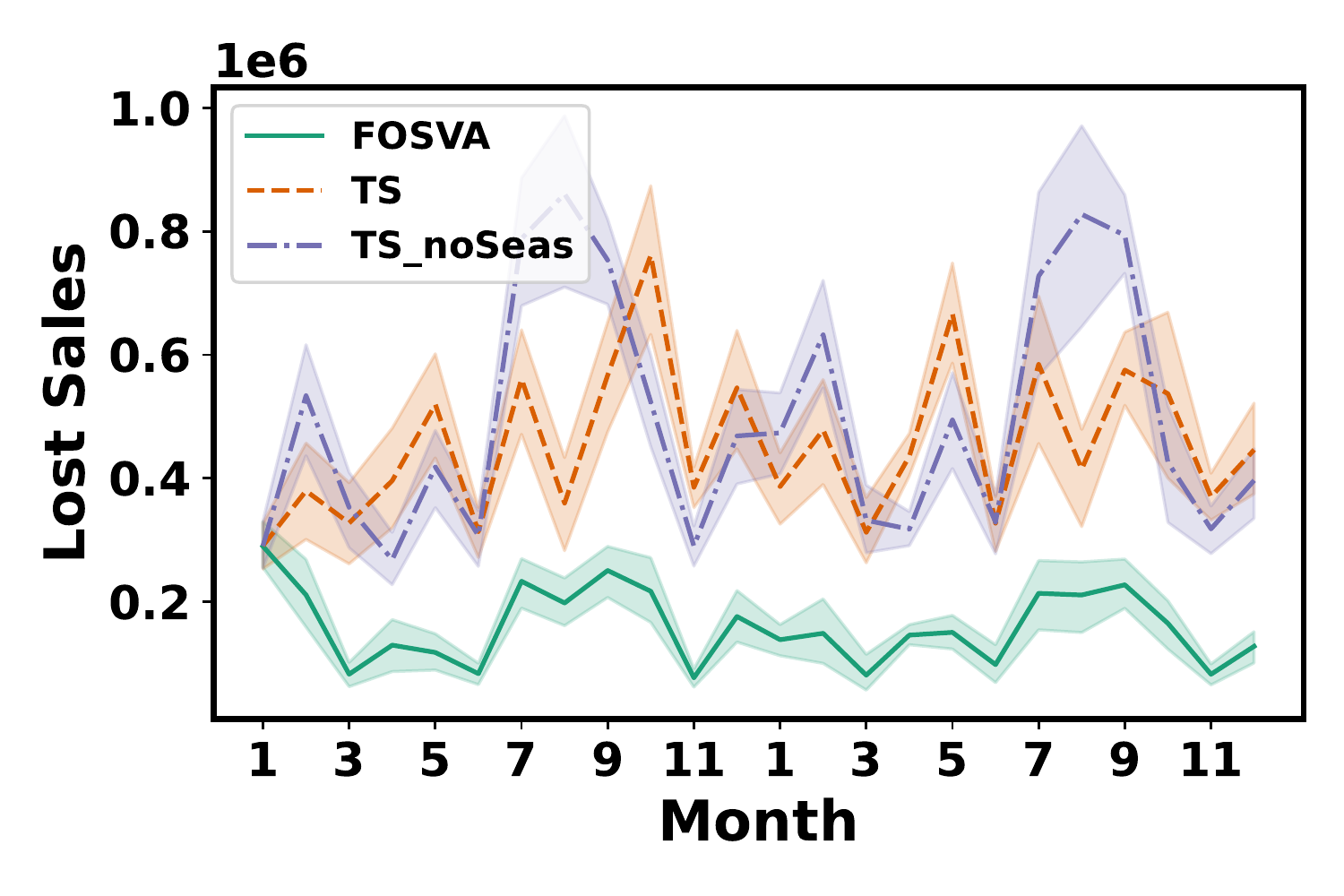}
\end{subfigure}

\caption{Mean and 95\% confidence intervals of production costs and lost sales for two-stage models over a 24-month horizon. Models exploit 10 years of historical data, with tightness $\gamma = 1.3$. Values are calculated over 10 independent out-of-sample replications of the rolling-horizon simulation.} \label{fig:TSstrategies}
\end{figure}

We start with strategies based on two-stage scenario trees (\TS, \TSn, 
and \FOSVA). Figure~\ref{fig:TSstrategies} shows the production costs (relative to the number of components produced by applying different rolling-horizon strategies) and lost sales during the 24-month out-of-sample simulation, together with a 95\% confidence interval calculated over 10 replications, assuming 
10 years of historical data are available and $\gamma =1.3$. The effect of initial inventory conditions disappears immediately, not affecting the periodicity of production volumes, which follow a periodical pattern as early as the second month. This result allows us to neglect the dependence on the initial state. As we will see, this phenomenon holds true for the other models and for the other experimental settings.

Regarding the characteristics of different two-stage scenario trees, in the case where we use observations without paying attention to seasonality (\TSn), production volumes consistently follow a non-seasonal strategy, leading to great lost sales. On the contrary, \TS properly follows seasonality but, since it cannot rely on a terminal value estimation, it has lower production volumes than \FOSVA. This results in the inability to cope with seasonal peaks, leading to lost sales. These observations remain true for the other configurations that are not visualized 
graphically (due to space constraints), but are reported in   Tables~\ref{tab:profits},~\ref{tab:inventory}, and \ref{tab:lostSales}.

In terms of mean profits (Table~\ref{tab:profits}), \FOSVA outperforms \TS and \TSn in all configurations. This is due to insufficient stocks of components, as can be seen in Table~\ref{tab:inventory}. In particular, \TS and \TSn have an average level of inventory lower than what is needed to meet demand. This stock shortage is consistently recorded in Table~\ref{tab:lostSales}, where we observe average lost sales that are often more than twice as high as the other methods (i.e., $>100\%$).

\begin{table}[tb]
\caption{Mean percentage profit with respect to the perfect information model, computed over 10 replications of the 24-month out-of-sample horizon. Values are displayed 
by amount of available data and capacity tightness.}\label{tab:profits}
\begin{center}
\begin{tabular}{c c | c c c c c c c c c}
&&\multicolumn{9}{c}{Policy}\\
Years & $\gamma$ & \FOSVA & \TS & \TSn & \MPd & \MPt & \MPq & \MSt & \MStt & \MStq \\ \hline
\multirow{4}{*}{3} & $1.3$ & 42.9& 17.5& 15.8& 40.6& 41.4& 40.6& 41.2& 41.9& 41.3\\
& $1.2$ & 43.1 & 17.2 & 16.1 & 40.5 & 41.2 & 40.5 & 40.8 & 41.5 & 40.5\\ 
& $1.1$ & 42.9& 16.7& 16.7& 40.3& 41.2& 40.8& 40.6& 41.1& 39.8\\ 
& $1.0$ &43.3& 16.5& 17.6& 41.3& 42.0& 41.8& 41.3& 41.2& 40.1\\
\hline
\multirow{4}{*}{5} & $1.3$ &48.4& 18.0& 16.5& 44.9& 46.8& 47.0& 45.1& 46.8& 46.5 \\
& $1.2$ &48.4& 18.0& 16.8& 44.6& 46.6& 46.9& 44.9& 46.7& 46.0\\ 
& $1.1$ &48.1& 17.8& 17.3& 44.9& 47.2& 47.1& 45.0& 46.8& 46.0\\ 
& $1.0$ &48.6& 18.0& 18.3& 46.0& 47.8& 47.7& 46.0& 47.2& 46.2\\
\hline
\multirow{4}{*}{10} & $1.3$ &49.6& 17.1& 15.5& 46.8& 48.4& 48.3& 47.8& 48.4& 47.8\\
& $1.2$ & 49.5& 17.1& 15.8& 46.8& 48.5& 48.6& 47.4& 48.5& 48.0\\ 
& $1.1$ & 49.6& 17.0& 16.3& 47.5& 50.1& 50.0& 48.1& 49.9& 49.1\\ 
& $1.0$ & 51.0& 17.0& 17.2& 50.1& 52.5& 52.4& 50.2& 51.9& 51.5\\
\hline
\end{tabular}
\end{center}
\end{table}

\begin{table}[tb]
\caption{Average inventory level with respect to the perfect information model, computed over 10 replications of the 24-month out-of-sample horizon. Values are displayed by capacity tightness and the amount of available data.}\label{tab:inventory}
\begin{center}
\begin{tabular}{c c | c c c c c c c c c}
&&\multicolumn{9}{c}{Policy}\\
Years & $\gamma$ & \FOSVA & \TS & \TSn & \MPd & \MPt & \MPq & \MSt & \MStt & \MStq \\ \hline
\multirow{4}{*}{3} & $1.3$ & 119.0 & 85.1 & 73.2 & 117.3 & 125.4 & 126.8 & 118.3 & 126.1 & 128.2 \\
& $1.2$ &117.6& 84.2& 72.8& 115.7& 123.6& 124.9& 116.2& 125.5& 128.3\\ 
& $1.1$ &115.6& 83.3& 72.7& 113.3& 123.9& 125.8& 113.8& 125.0& 128.2\\ 
& $1.0$ &113.7& 83.3& 73.4& 111.7& 122.1& 123.6& 111.7& 122.5& 125.8\\
\hline
\multirow{4}{*}{5} & $1.3$ &120.5& 83.3& 73.5& 120.9& 131.8& 134.0& 123.1& 131.9& 135.4 \\
& $1.2$ &118.9& 82.7& 73.3& 118.7& 130.0& 132.0& 120.6& 130.4& 134.3\\ 
& $1.1$ &116.2& 81.8& 73.1& 116.0& 129.2& 131.7& 117.2& 129.3& 133.6\\ 
& $1.0$ &113.3& 81.9& 73.7& 114.5& 126.2& 128.2& 114.4& 125.9& 129.8\\
\hline
\multirow{4}{*}{10} & $1.3$ &120.2& 79.6& 72.8& 120.2& 133.6& 135.5& 123.2& 133.9& 138.0 \\
& $1.2$ &117.6& 78.4& 71.9& 117.2& 130.3& 132.4& 119.9& 130.8& 135.5\\ 
& $1.1$ &114.5& 77.4& 71.4& 114.2& 129.2& 131.9& 116.2& 128.7& 134.2\\ 
& $1.0$ &111.2& 76.9& 71.5& 111.9& 124.9& 127.3& 112.6& 124.4& 129.4\\
\hline
\end{tabular}
\end{center}
\end{table}

\begin{table}[tb] 
\caption{Lost sales percentage with respect to all the models, computed over 10 replications of the 24-month out-of-sample horizon. Values are displayed by capacity tightness and the amount of available data.}\label{tab:lostSales}
\begin{center}
\begin{tabular}{c c | c c c c c c c c c}
&&\multicolumn{9}{c}{Policy}\\
Years & $\gamma$ & \FOSVA & \TS & \TSn & \MPd & \MPt & \MPq & \MSt & \MStt & \MStq \\ \hline
\multirow{4}{*}{3} & $1.3$ & -13.6& 101.2& 124.0& -9.1& -22.9& -24.3& -11.2& -24.4& -26.7 \\
& $1.2$ & -13.8& 96.7& 117.1& -8.6& -21.2& -22.4& -9.5& -24.0& -26.5\\ 
& $1.1$ & -11.9& 90.7& 107.3& -5.9& -21.5& -23.4& -6.7& -23.3& -25.6\\ 
& $1.0$ & -10.4& 77.4& 89.7& -5.2& -18.6& -20.0& -5.2& -19.4& -22.0\\
\hline
\multirow{4}{*}{5} & $1.3$ & -11.2& 127.5& 155.8& -12.0& -35.3& -39.4& -16.8& -35.3& -41.1 \\
& $1.2$ & -10.9& 118.8& 144.9& -10.2& -32.5& -36.4& -14.3& -33.3& -39.2\\ 
& $1.1$ & -8.0& 108.5& 130.7& -7.8& -31.4& -35.0& -10.2& -31.5& -37.1\\ 
& $1.0$ & -6.2& 89.5& 107.1& -7.5& -26.0& -28.7& -8.4& -25.5& -29.9\\
\hline
\multirow{4}{*}{10} & $1.3$ &-12.5& 148.5& 170.7& -12.4& -42.6& -45.7& -20.7& -43.0& -49.4\\
& $1.2$ &-11.9& 139.9& 160.4& -11.2& -39.6& -43.1& -18.0& -40.3& -47.3\\ 
& $1.1$ & -9.2& 129.7& 146.9& -8.4& -38.6& -41.9& -14.2& -37.8& -45.8\\ 
& $1.0$ & -7.7& 109.4& 122.4& -9.1& -32.1& -35.2& -12.0& -31.4& -38.1\\
\hline
\end{tabular}
\end{center}
\end{table}

\subsection{Multistage models}

\begin{figure}[t]
\centering
\begin{subfigure}[]{0.49\textwidth}
\includegraphics[width=\textwidth]{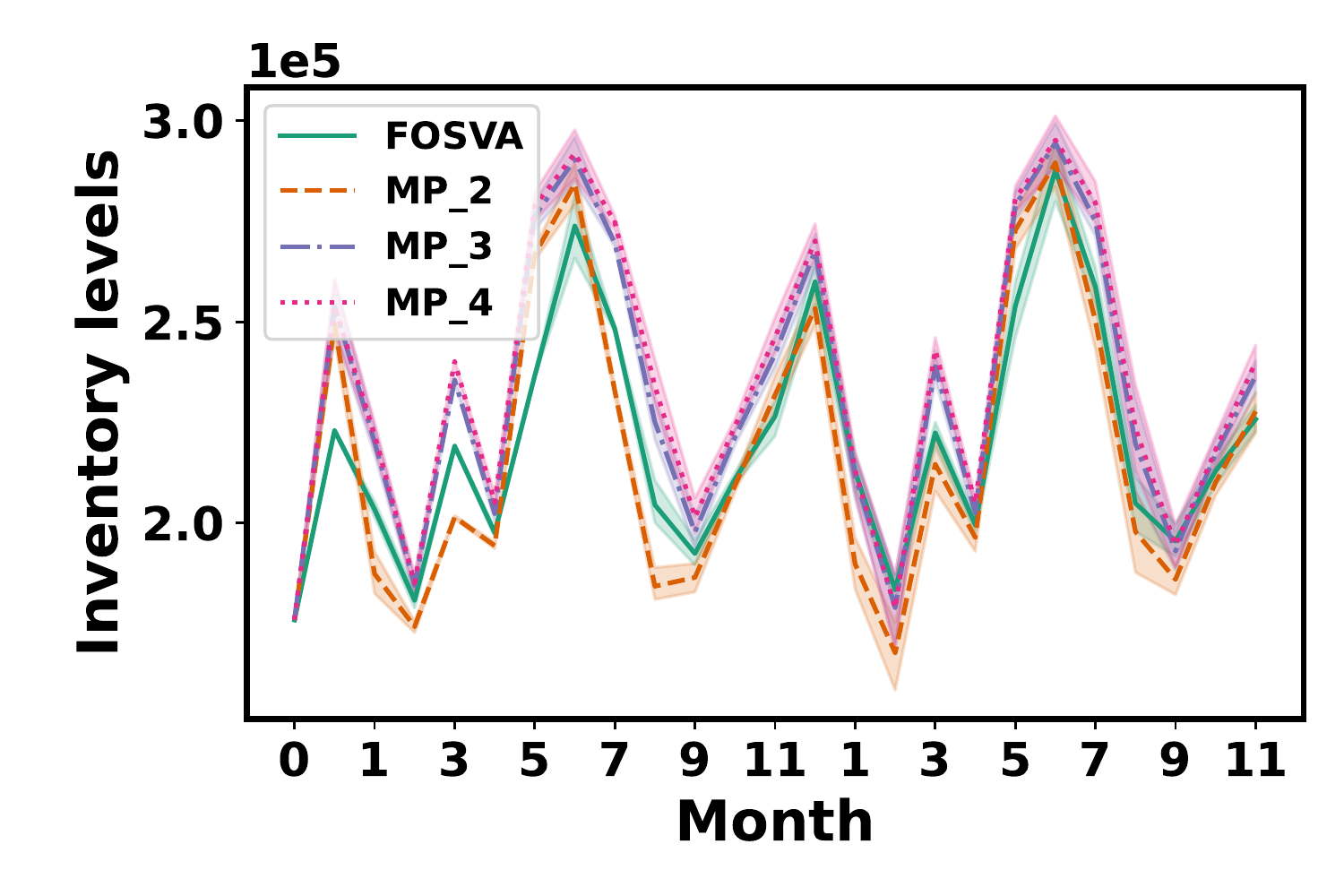}
\end{subfigure}
\,
\begin{subfigure}[]{0.49\textwidth}
\includegraphics[width=\textwidth]{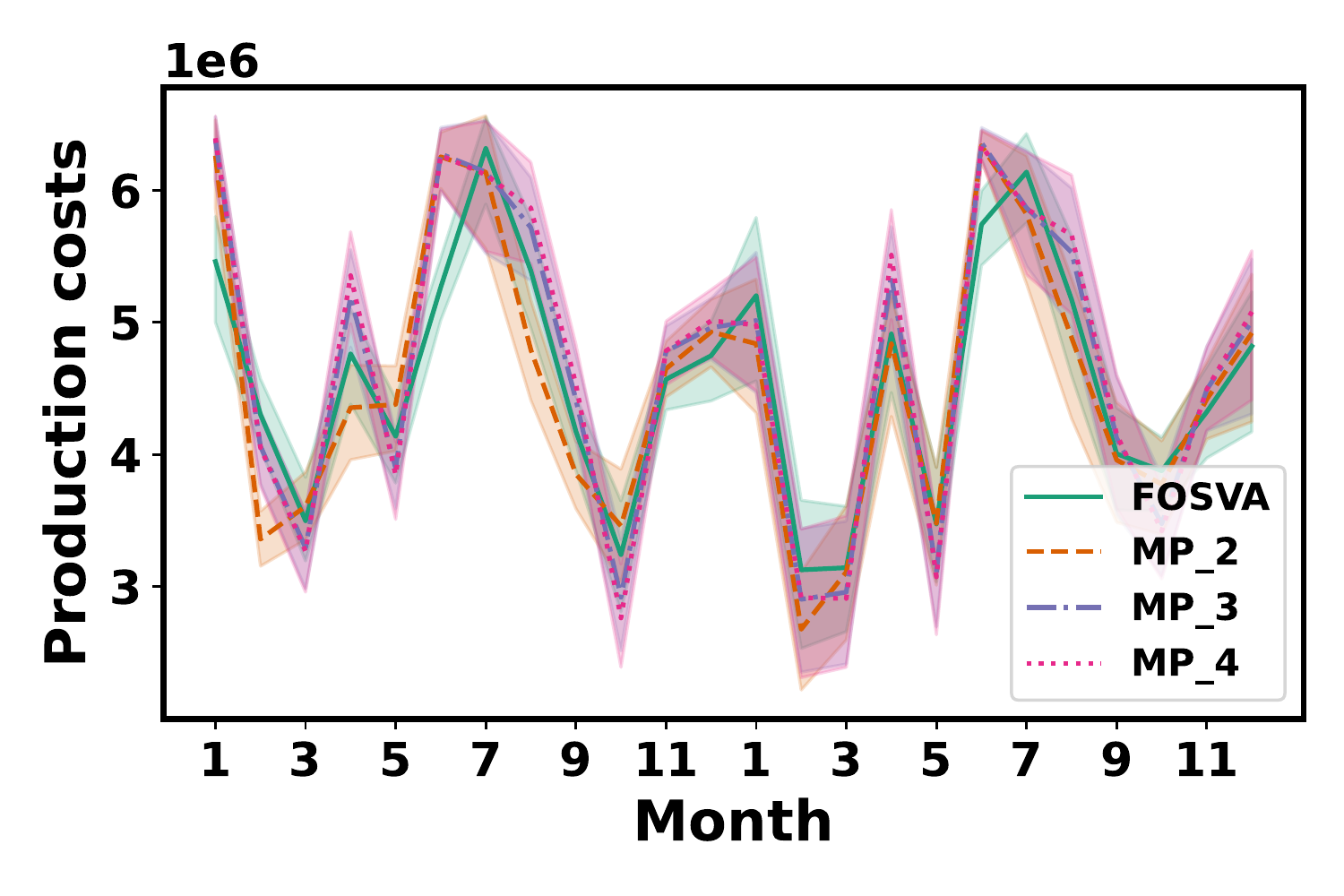}
\end{subfigure}
\\
\begin{subfigure}[]{0.49\textwidth}
\includegraphics[width=\textwidth]{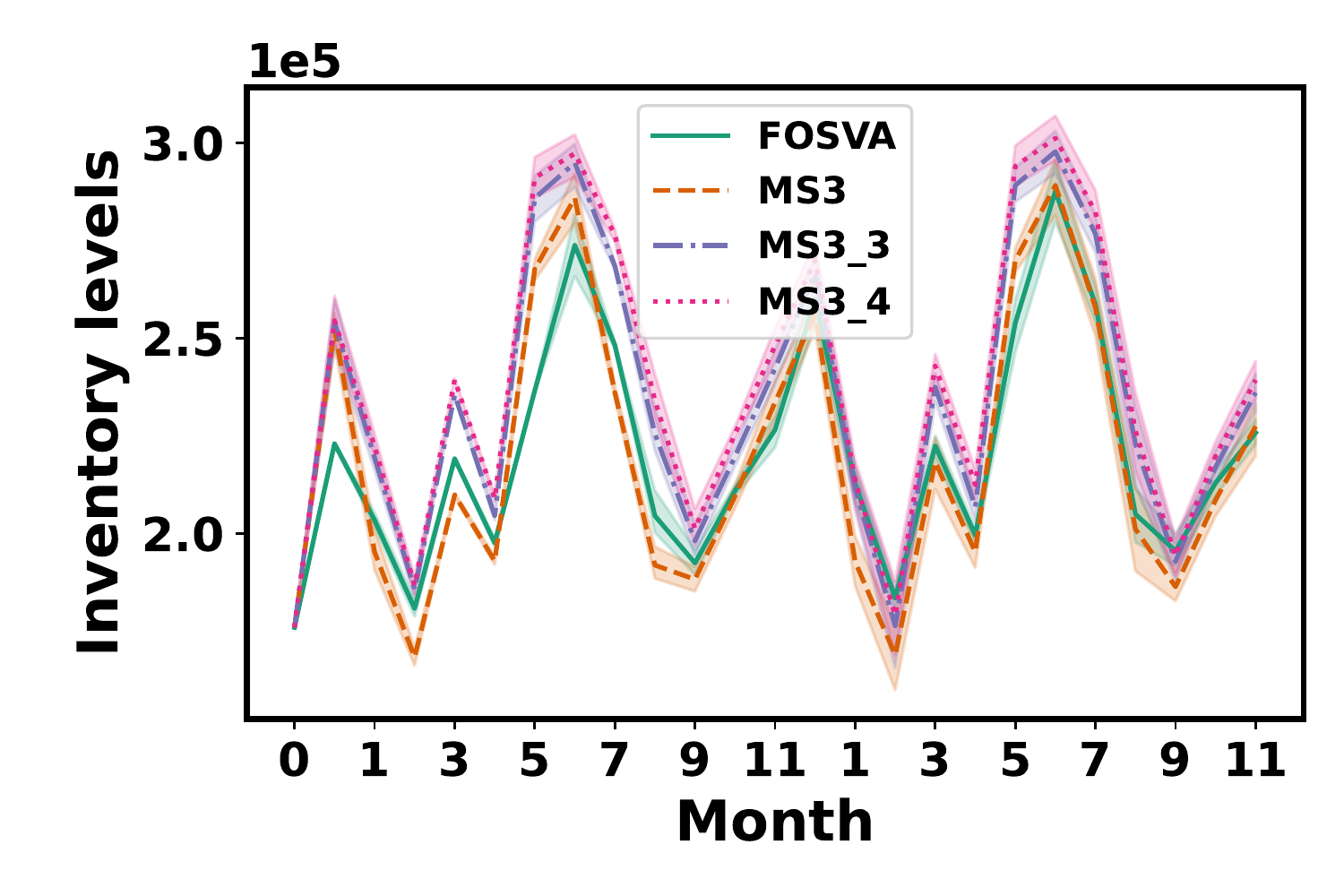}
\end{subfigure}
\,
\begin{subfigure}[]{0.49\textwidth}
\includegraphics[width=\textwidth]{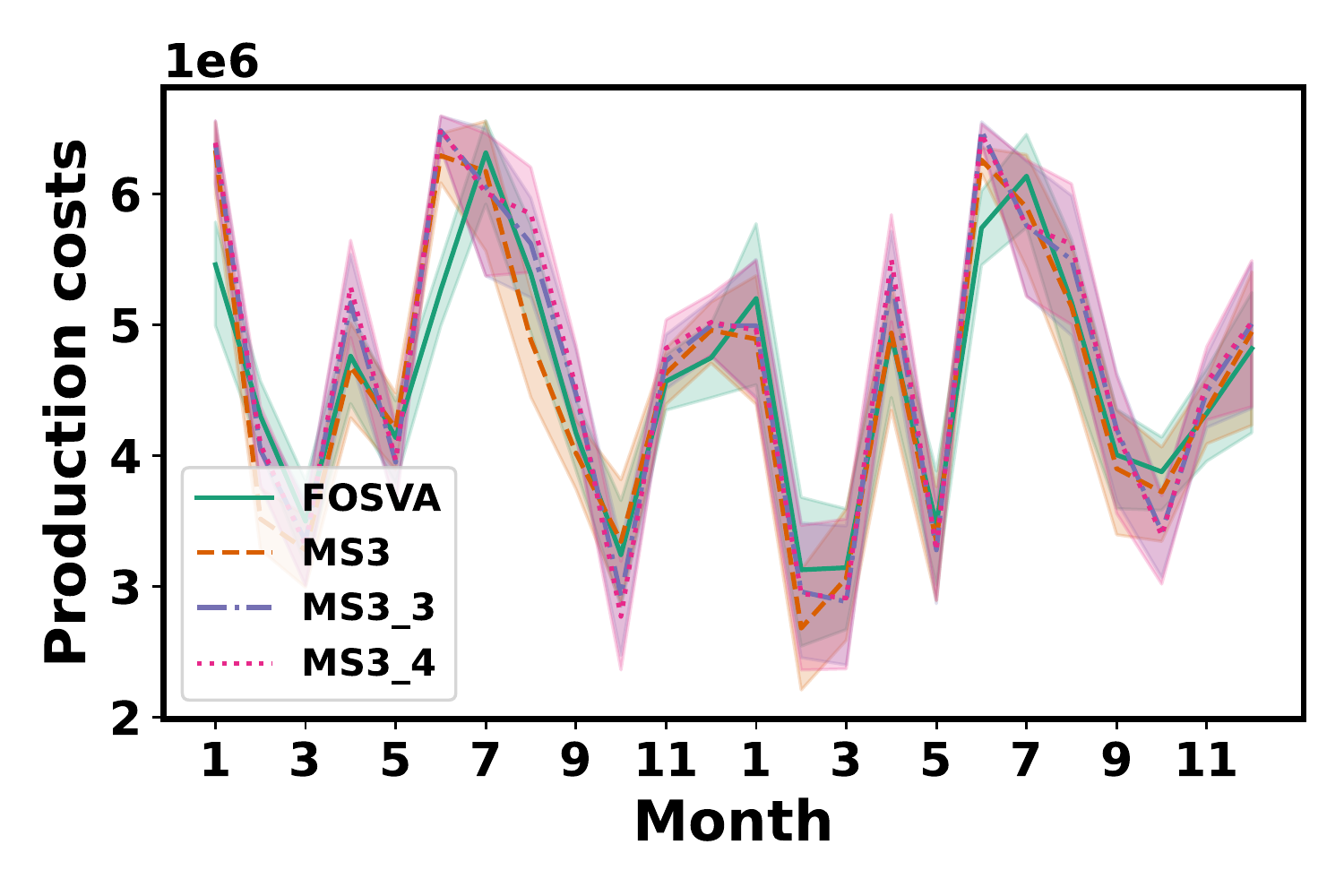}
\end{subfigure}

\caption{Mean and 95\% confidence intervals of production costs and inventory levels for multistage models over a 24-month horizon. Models are based on 3 years of historical data and a tightness $\gamma = 1.3$. Values are calculated over 10 independent out-of-sample replications of the rolling-horizon simulation.} \label{fig:startEff}
\end{figure}

Similarly to cases featuring a two-stage scenario tree, 
in order to check whether the results discussed are independent of the 
initial inventory conditions, it is useful to investigate when the production volume becomes periodic. Figure~\ref{fig:startEff} shows the mean of production costs and component quantities in the inventory for 10 replications of the 24-month out-of-sample horizon, with $\gamma=1.3$, when 3 years of historical demand are available. 
The 95\% confidence intervals are reported. 
The impact of initial inventory conditions fades away after a short transient phase, resulting in periodic production patterns. This behavior can also be observed in configurations that we do not report here 
due to space limitations.

Considering the profits of the different methods in various configurations (Table~\ref{tab:profits}), we observe that \FOSVA has superior performance in almost all configurations except for low tightness ($\gamma = 1.0, 1.1)$, and a large amount 
of available historical observations (10 years). In general, \MStt and \MPt achieve the second-best performance, being 1\% - 2\% worse than \FOSVA. It is worth noting that these results come from different production strategies. Using Figure~\ref{fig:startEff} as a graphical illustration, we note that \FOSVA generally maintains a lower level of inventory than \MStt and \MPt, approaching those of the methods with a two-periods lookahead (\MSt and \MPd). From Table~\ref{tab:inventory}, we can assert that a similar pattern holds for all configurations. Furthermore, we see that as the horizon of the tree increases, the inventory levels increase accordingly. Nevertheless, a longer view and higher inventory are not always favorable; in fact, \MStq and \MPq hold an excessive inventory, ranking behind \MStt and \MPt in terms of profit, due to the associated holding costs.
Lost sale costs are reported in Table~\ref{tab:lostSales}. They follow the same pattern of the average inventory levels, decreasing for all models (i.e., more negative percentages) as tightness increases. However, a good balance between the number of components in inventory and lost sales is not achieved for too long horizons, resulting in lower profits.


\subsection{Deterministic model with safety stocks}\label{subsec:ss}

Table~\ref{tab:SSprofits} reports the results of the safety-stocks-based heuristics for $\alpha = 0, 10, 15$, and $50$. We determine the empirical distribution from which $\alpha$ is computed by aggregating the observations without considering seasonality. This approach is beneficial because the number of observed scenarios is too limited, entailing unreasonably large safety stock values, driven by an insufficiently sized data pool. Consequently, the number of observed years does not produce a significant effect, instead, it slightly modifies both the performance and inventory level.

Notably, \SSq{0} does not account for a safety stock value, but rather it embodies a model that acknowledges uncertainty solely through the scenarios mean.
As $\alpha$ increases, the inventory level correspondingly rises, reaching a similar value to \FOSVA when  $\alpha = 10$. With this value, the best results are achieved. All of the methods but \SSq{50} outperform the conventional two-stage methods (\TS and \TSn), reaching an average percentage profit (with respect to full information) of 31.17\%, against 17.32\% by \TS and 16.65\% by \TSn, respectively. This difference is primarily due to the combination of safety stocks and a lookahead strategy provided by the model. In fact, despite its deterministic approach, it offers a significant advantage when a minimum stock level is incorporated. Nevertheless, the best safety stock strategy (\SSq{10}) achieves results worse than all the other methods (i.e., \FOSVA, \MPd, \MPt, \MPq, \MSt, \MStt, \MStq), proving that safety stock models under-perform problems embedding multiple future scenarios, 
which provide us with a better hedge against demand uncertainty. Lastly, in contrast to other strategies, safety stock-based heuristics tend to decline in performance as capacity tightness diminishes. This trend is contradicted only in settings with  $\alpha = 50$ and 10 years of available observations. Although the safety stock required by the heuristic remains constant relative to the tightness, the flexibility of the model is hampered when the capacity is smaller.

\begin{table}[tb]
\caption{Mean percentage profit and average inventory level of the safety stock policies (\SSq{quantile}) with respect to the perfect information model, computed over 10 replications of the 24-month out-of-sample horizon. Values are displayed by amount of available data and capacity tightness.}\label{tab:SSprofits}
\begin{center}
\begin{tabular}{c c | c c c c c c c c}
&&\multicolumn{8}{c}{Policy}\\
\hline
&&\multicolumn{4}{c|}{Profit}&\multicolumn{4}{c}{Inventory}\\
Years & $\gamma$ & \SSq{0} & \SSq{10} & \SSq{25} & \multicolumn{1}{c|}{\SSq{50}} & \SSq{0} & \SSq{10} & \SSq{25} & \SSq{50} \\ \hline
\multirow{4}{*}{3} & $1.3$ & 32.0& 35.3& 32.4& \multicolumn{1}{c|}{16.0}& 106.8& 116.7& 146.2& 181.5\\
& $1.2$ & 31.9& 34.9& 31.6& \multicolumn{1}{c|}{15.1}& 108.7& 118.3& 147.2& 181.0\\ 
& $1.1$ & 30.2& 33.3& 30.1& \multicolumn{1}{c|}{13.6}& 113.6& 122.1& 147.4& 178.8\\ 
& $1.0$ &25.2& 29.0& 27.6& \multicolumn{1}{c|}{13.0}& 122.6& 129.3& 150.7& 177.0\\
\hline
\multirow{4}{*}{5} & $1.3$ &32.2& 35.9& 32.5& \multicolumn{1}{c|}{15.0}& 108.0& 118.7& 148.8& 184.4 \\
& $1.2$ &32.0& 35.5& 31.6& \multicolumn{1}{c|}{14.5}& 110.5& 121.0& 150.1& 184.0\\ 
& $1.1$ &29.4& 33.3& 30.5& \multicolumn{1}{c|}{14.2}& 116.7& 125.7& 151.4& 182.4\\ 
& $1.0$ &23.5& 28.4& 28.3& \multicolumn{1}{c|}{12.6}& 126.6& 133.2& 153.7& 180.1\\
\hline
\multirow{4}{*}{10} & $1.3$ &31.0& 35.7& 32.4& \multicolumn{1}{c|}{16.0}& 106.8& 118.4& 147.5& 184.6\\
& $1.2$ & 30.6& 35.4& 32.0& \multicolumn{1}{c|}{16.0}& 108.6& 119.9& 147.8& 183.2\\ 
& $1.1$ & 28.3& 33.8& 32.1& \multicolumn{1}{c|}{16.3}& 113.5& 122.9& 147.6& 180.4\\ 
& $1.0$ & 23.6& 30.1& 30.8& \multicolumn{1}{c|}{17.1}& 121.5& 128.7& 149.1& 175.4\\
\hline
\end{tabular}
\end{center}
\end{table}
%

\section{Conclusions and directions for further research}\label{sec:conclusions}

In this paper, we have presented 
different stochastic programming models for solving finite-capacity ATO production planning 
problems. Besides proposing an open-source, comprehensive framework for the simulation and optimization of these types of problems, we have focused on the behavior of different models and on how they use the available information. Moreover, we have presented 
a methodology for assigning a value to end-of-horizon stocks and reducing the myopia of (two-stage) models with short horizons. Furthermore, we have investigated the effects of correlation and seasonality, defining hierarchical structures of product families, which are a natural occurrence in an ATO setting. 

The addition of a terminal value improved performance considerably, 
leading to solutions of better quality than those obtained by models based on a longer planning horizon. This result is particularly interesting, considering that the approximation of the 
value of the terminal state is based on a decomposition with respect to components and neglects the interaction among components through the bills of materials. Moreover, we have shown that considering a longer planning horizon does not necessarily improve performance. 
Similar conclusions are reported, in a different setting, by \citep{BirgeTree}. 

From a practical viewpoint, the proposed approach could be useful to provide support for two-level MPS/FAS procedures in MRP systems.
From a more general algorithmic viewpoint, 
it can be interpreted as an approximate dynamic programming strategy, integrating a limited lookahed tree with an approximation of the state value function.
Indeed, the results obtained by adopting a concave, piecewise linear
approximation of the value function of the terminal value seem promising and applicable to other contexts. Future studies will consider the impact of such a strategy in applications like the classical capacitated lot-sizing problem.

\FloatBarrier


\end{document}